\RequirePackage{ifpdf}
\ifpdf 
\documentclass[pdftex]{sigma}
\else
\documentclass{sigma}
\fi

\usepackage{mathrsfs}
\usepackage[mathscr]{euscript}

\newtheorem{thm}{Theorem}[section]
\newtheorem{lem}[thm]{Lemma}
\newtheorem{prop}[thm]{Proposition}
\newtheorem{cor}[thm]{Corollary}

\theoremstyle{definition}
\newtheorem{defn}[thm]{Definition}
\newtheorem{ex}[thm]{Example}
\newtheorem{rem}[thm]{Remark}

\newcommand{\osp}{\mathfrak{osp}}
\newcommand{\spk}{\mathfrak{sp}}

\newcommand{\ok}{\mathfrak{o}}

\newcommand{\Sk}{\mathfrak{S}}
\newcommand{\g}{\mathfrak{g}}

\newcommand{\hk}{\mathfrak{h}}

\newcommand\Bc{\mathcal{B}}
\newcommand\Lc{\mathcal L}

\newcommand\Ac{\mathcal A}

\newcommand\Cc{\mathcal C}
\newcommand\Mcc{\mathcal M}
\newcommand\Sc{\mathcal S}

\newcommand\Uc{\mathcal U}
\newcommand\Vc{\mathcal V}
\newcommand\Zc{\mathcal Z}

\newcommand{\Ab}{\mathsf{A}}
\newcommand{\Bb}{\mathsf{B}}
\newcommand{\Db}{\mathsf{D}}
\newcommand{\Wb}{\mathsf{W}}
\newcommand{\Sb}{\mathsf{S}}
\newcommand{\Mb}{\mathsf{M}}
\newcommand{\Zb}{\mathsf{Z}}

\newcommand\CC{\mathbb C}
\newcommand\NN{\mathbb N}
\newcommand\ZZ{\mathbb Z}

\newcommand{\Tr}{\operatorname{Tr}}

\newcommand{\SP}{\operatorname{SP}}

\newcommand{\Diff}{\operatorname{Dif\/f}}
\newcommand{\ad}{\operatorname{ad}}

\newcommand{\expo}{\operatorname{exp}}

\renewcommand\dfrac{\displaystyle \frac}

\newcommand{\sym}{\operatorname{Sym}}
\newcommand{\Id}{\operatorname{Id}}
\renewcommand\hat\widehat
\renewcommand\tilde\widetilde 
\newcommand{\spa}{\operatorname{span}}

\newcommand{\Str}{\operatorname{Str}}

\newcommand\sta{{ \ \star \ }}

\newcommand\tstu{{\mathop{\star}\limits_1 \, }}
\newcommand\tsti{{\mathop{\star}\limits_t}}

\newcommand\Wedge{\bigwedge}

\newcommand{\Cs}{\mathscr C}
\newcommand{\Es}{\mathscr E}

\newcommand{\Ps}{\mathscr P}

\newcommand{\apl}{{\mathcal{A}_\lambda}}
\newcommand{\apln}{{\mathcal{A}_\lambda(n)}}
\newcommand{\agl}{\mathcal{A}_\Lambda}
\newcommand{\Agl}{A_\Lambda}
\newcommand{\agln}{\mathcal{A}_\Lambda(n)}
\newcommand{\azn}{\mathcal{A}_0(n)}

\newcommand{\ut}{\mathcal{U}_\vartheta}
\newcommand{\dzd}{\operatorname{deg}_{\ZZ_2}}
\newcommand{\zero}{\overline{0}}
\newcommand{\um}{\overline{1}}
\newcommand{\bl}{{\scriptscriptstyle{\Lc}}}
\newcommand{\ze}{{\scriptscriptstyle{\overline{0}}}}
\newcommand{\un}{{\scriptscriptstyle{\overline{1}}}}
\newcommand{\ib}{{\scriptscriptstyle{\overline{i}}}}
\newcommand{\jb}{{\scriptscriptstyle{\overline{j}}}}

\newcommand{\gO}{\g_{\ze}}
\newcommand{\gI}{\g_{\un}}

\newcommand{\zdtimes}{\mathop{\otimes}\limits_{\ZZ_2}}

\newcommand{\SM}{\mathsf{SM}}

\newcommand{\lw}{{\scriptscriptstyle{\Wedge}}}
\newcommand{\lsb}{{\scriptscriptstyle{\Sb}}}

\begin{document}

\allowdisplaybreaks

\renewcommand{\PaperNumber}{028}

\renewcommand{\thefootnote}{$\star$}

\FirstPageHeading

\ShortArticleName{Hochschild Cohomology and Deformations of
  Clif\/ford--Weyl Algebras}

\ArticleName{Hochschild Cohomology and Deformations\\ of
  Clif\/ford--Weyl Algebras\footnote{This paper
is a contribution to the Special Issue on Deformation
Quantization. The full collection is available at
\href{http://www.emis.de/journals/SIGMA/Deformation_Quantization.html}{http://www.emis.de/journals/SIGMA/Deformation\_{}Quantization.html}}}

\Author{Ian M. MUSSON~$^\dag$, Georges PINCZON~$^\ddag$ and Rosane USHIROBIRA~$^\ddag$}

\AuthorNameForHeading{I.M. Musson, G. Pinczon and R. Ushirobira}

\Address{$^\dag$~Department of Mathematical Sciences, University of
  Wisconsin-Milwaukee,\\
  \hphantom{$^\dag$}~Milwaukee, WI 53201-0413,  USA}

\EmailD{\href{mailto:musson@uwm.edu}{musson@uwm.edu}}
\URLaddressD{\url{http://www.uwm.edu/~musson/}}

\Address{$^\ddag$~Institut de Math\'ematiques de Bourgogne,
  Universit\'e de Bourgogne,\\
   \hphantom{$^\ddag$}~B.P.\ 47870, F-21078 Dijon Cedex, France}

\EmailD{\href{mailto:Georges.Pinczon@u-bourgogne.fr}{Georges.Pinczon@u-bourgogne.fr},
  \href{mailto:Rosane.Ushirobira@u-bourgogne.fr}{Rosane.Ushirobira@u-bourgogne.fr}}
\URLaddressD{\url{http://www.u-bourgogne.fr/monge/phy.math/members/pinczon.htm},\\
\hspace*{13.5mm}\url{http://www.u-bourgogne.fr/rosane.ushirobira}}

\ArticleDates{Received October 01, 2008, in f\/inal form February 25,
2009; Published online March 07, 2009}

\Abstract{We give a complete study of the Clif\/ford--Weyl algebra
  ${\mathcal C}(n,2k)$ from Bose--Fermi statistics, including Hochschild
  cohomology (with coef\/f\/icients in itself). We show that ${\mathcal C}(n,2k)$
  is rigid when $n$ is even or when $k \neq 1$. We f\/ind all
  non-trivial deformations of ${\mathcal C}(2n+1,2)$ and study their
  representations.}


\Classification{16E40; 16G99; 16S80; 17B56; 17B10; 53D55}

\pdfbookmark[1]{Introduction}{Introduction}
\section*{Introduction}\label{section0}

Throughout the paper, the base f\/ield is $\CC$. As usual in
superalgebra theory, we denote the ring $\ZZ/ 2 \ZZ$ by $\ZZ_2$.

Let $\Cc(n)$ be the Clif\/ford algebra with $n$ generators and
$\Wb_{2k}$ be the Weyl algebra with $2k$ generators. Denote by $V_\ze$
the vector space spanned by the generators of~$\Cc(n)$. Elements of~$V_\ze$ will be called {\em Fermi-type operators}. Similarly, let
$V_\un$ be the vector space spanned by the generators of $\Wb_{2k}$.
Elements of $V_\un$ will be called {\em Bose-type operators}.

There exist $\ZZ_2$-gradations on $\Cc(n)$ and $\Wb_{2k}$ such that
Fermi and Bose-type operators all have degree one. The {\em
  Clif\/ford--Weyl algebra} is:
\[\Cc(n, 2k) := \Cc(n) \otimes_{\ZZ_2} \Wb_{2k},\]
where $\otimes_{\ZZ_2}$ is relative to these gradations. It unif\/ies
Fermi and Bose-type operators in a unique algebra: as elements of
$\Cc(n,2k)$, they anti-commute. There is a $\ZZ_2$-gradation on
$\Cc(n,2k)$ extending the natural gradation of $V=V_\ze \oplus V_\un$,
and a corresponding structure of Lie superalgebra. Palev has shown
that $V$ generates a sub-superalgebra of $\Cc(n,2k)$ isomorphic to
$\osp(n+1,2k)$, and introduced corresponding parastatistics relations
\cite{Palev}. This was an outcome of previous results by Wigner
\cite{Wi}, Green \cite{Gre} and others (see \cite{FF}). It gives an
elegant algebraic interpretation of the Green ansatz, using the Hopf
structure of the enveloping algebra of $\osp(n+1,2k)$ \cite{Palev94},
and introduces a construction of representations of parastatistics
relations by Verma modules of~$\osp(n+1,2k)$. It also gives an idea of
what deformed (quantum) parastatistics could be: replace $\Cc(n,2k)$
by a ``quantum deformation'', which still has a Hopf structure. This
idea was developed by Palev himself~\cite{Palev94} and other authors.

The purpose of the present paper is to return to the f\/irst steps of
these theories: Clif\/ford--Weyl algebras. At this level, there is a very
natural question: does there exist non-trivial deformations of
$\Cc(n,2k)$? By a deformation, we mean a formal one, in the sense of
Gerstenhaber theory \cite{G64}. It is well-known that the answer is no
for $\Cc(n,0) = \Cc(n)$ and $\Cc(0,2k) = \Wb_{2k}$, but nothing was
done in the general case. We shall answer the question, but this is
not our only goal. We also want to introduce Clif\/ford--Weyl algebras in
a deformation quantization framework, emphasize their periodicity
behavior and how it can be used, explain where Palev's theorem comes
from, and so on.

Before describing the content of the paper, let us answer the initial
question: $\Cc(2n,2k)$ is rigid, for all $n$, $k$, $\Cc(2n+1,2k)$ is
rigid if, and only if, $k \neq 1$, so the answer is no in these
cases. In the case of $\Cc(2n+1,2)$ there exist non-trivial
deformations, that we completely describe in the paper, including
their representations.

Let us give some details of our main results. In Section
\ref{Section00}, we recall well-known properties of Clif\/ford and Weyl
algebras needed in the paper. In particular, we recall the
deformation quantization construction of the Weyl algebra
(resp.\ Clif\/ford algebra) through the Moyal product (resp.\ a Moyal-type
product).

In Section~\ref{Section01}, we construct the Clif\/ford--Weyl algebra $\Cc(n,2k)$ by a
similar deformation quantization procedure, as a deformation of the
super-exterior algebra of the $\ZZ_2$-graded vector space $V = V_\ze
\oplus V_\un$ with an explicit Moyal-type formula for the
$\sta$-product. From this construction, $\Cc(n,2k)$ is a $\ZZ_2 \times
\ZZ_2$-graded algebra, with natural left and right $\ZZ_2$-gradations.

We show in Section~\ref{section3}, that Clif\/ford--Weyl algebras have a periodic
behavior, very similar to the well-known, and useful, periodicity of
Clif\/ford algebras:

\medskip

\noindent
{\bf Periodicity Lemma 1.}
\[
\Cc(2m+n, 2k) \simeq \Cc(2m) \otimes \Cc(n,2k).
\]

This simple Lemma has many consequences. For a given $r \in \NN^*$,
denote by $\Mcc_r(A)$ the $r \times r$-matrix algebra with
coef\/f\/icients in an algebra $A$. One has:

\medskip

\noindent
{\bf Theorem 1.}
\begin{enumerate}\itemsep=0pt

\item $\Cc(2n,2k) \simeq \Mcc_{2^n} \left( \Wb_{2k} \right)$.

\item $\Cc(2n+1,2k) \simeq \Mcc_{2^n} \left( \Cc(1,2k) \right)$.

\end{enumerate}

Note that $\Cc(1,2k)$ is the algebra generated by $\Wb_{2k}$ and the
parity operator of the metaplectic (oscillator) representation. In
other words, $\Cc(1,2k)$ is the smash product $\Sk_2 \ltimes \Wb_{2k}$
where $\Sk_2$ is the group $\{-1, 1\}$ (see Remark~\ref{3.4}). By
Periodicity Lemma~1, it results that, though $\Cc(2n, 2k)$ has been
def\/ined by anti-commuting Fermi and Bose-type operators, sitting in
$\osp(2n+1,2k)$, it can also be def\/ined by a new set of commuting
generators, also of Fermi and Bose-type, no longer sitting in
$\osp(2n+1, 2k)$, but generating a superalgebra of type $\ok(2n+1)
\times \osp(1,2k)$. This suggest that the enveloping algebra
$\mathcal{U}(\mathfrak{o}(2n+1) \times
\mathfrak{osp}(1,2k))$ could be used for
parastatistics. We shall not go further into this point in the present
paper.

A second consequence, and the key for our purpose (to deform
Clif\/ford--Weyl algebras) is~the~natural Morita equivalence (see
\cite{Loday}) between $\Wb_{2k}$ and $\Cc(2n,2k)$ (resp.\ $\Cc(1,2k)$
and $\Cc(2n+1, 2k)$). We shall see some consequences in Sections~\ref{section5}--\ref{section8}.

A new proof of Palev's theorem is given in Section~\ref{section4}:

\medskip

\noindent
{\bf Theorem 2 (\cite{Palev}).} {\it The sub-superalgebra of $\Cc(n,2k)$ generated by Fermi and
  Bose-type operators is isomorphic to $\osp(n+1,2k)$.}

\medskip

The original proof uses a direct comparison of commutation rules, and
does not really explain why the result exists. This is what we want to
do, and the reason for giving this new proof. We explain it in a few
words. First, we remark that $\Cc(n,2k)$ is $\ZZ_2 \times
\ZZ_2$-graded, and that the superbracket used in the Theorem 2 is the
one def\/ined by the right $\ZZ_2$-gradation. There exists a twisted
adjoint action $\ad'$ of $\Cc(n,2k)$ on itself, coming from the left
$\ZZ_2$-gradation, and a~supersymmetric bilinear form of type
$(n+1,2k)$ on $H:= \CC \oplus V$ (graded by $H_\ze = \CC \oplus
V_\ze$, $H_\un = V_\un$) coming from the natural super-Poisson
bracket. Then $\g = V \oplus [V,V]$ is a sub-superalgebra of
$\Cc(n,2k)$, $H$ is $\ad'(\g)$-stable, the bilinear form is
$\ad'(\g)$-invariant, and the result follows.

In Section~\ref{section5}, we go back to the initial question: to deform
$\Cc(n,2k)$. As well-known in deformation theory, the f\/irst step is to
study the Hochschild cohomology of $\Cc(n,2k)$ with coef\/f\/icients in
itself (see Appendix~\ref{appendix1} for generalities in Hochschild cohomology and
its relations with deformation theory). Since Morita equivalent
algebras have isomorphic Hochschild cohomology~\cite{Loday}, using the
Periodicity Lemma~1, it results that:

\medskip

\noindent
{\bf Theorem 3.}
\begin{enumerate}\itemsep=0pt
\item[$1.$] $H^\ell(\Cc(2n,2k)) =  \{0\}$ \textit{if} $\ell  > 0$.

\item[$2.$] $H^\ell(\Cc(2n+1,2k)) = H^\ell(\Cc(1,2k))$, \textit{for all} $\ell$.
\end{enumerate}

So we have a partial answer to our question: $\Cc(2n,2k)$ cannot be
non-trivially deformed, and we are left with the case of
$\Cc(1,2k)$. As we mentioned before, $\Cc(1,2k)$ is a smash product
$\Sk_2 \ltimes \Wb_{2k}$, but the cohomology of $\Sk_2 \ltimes
\Wb_{2k}$ is known, as a particular case of general results in~\cite{Alev} and \cite{P07}, that give the cohomology of $G \ltimes
\Wb_{2k}$, when $G$ is a f\/inite subgroup of $\SP(2k)$. We obtain:

\medskip

\noindent
{\bf Theorem 4.}
\begin{enumerate}\itemsep=0pt

\item[$1.$] {\it If $\ell >0$ and $\ell \neq 2k$, then}
  \[ H^\ell (\Cc(2n+1,2k)) = \{0\}.\]

\item[$2.$] {\it $\dim\left(H^{2k}(\Cc(2n+1,2k))\right) = 1.$ Denoting by
  $\omega_1, \dots, \omega_{2k+1}$ the basic Fermi-type operators and
  by $s$ the canonical symplectic form on $V_\un$, then there exists a
  $2k$-cocycle $\Omega$ such that $H^{2k} (\Cc(2n+1,2k)) = \CC
  \Omega$,
\begin{gather*}
   \Omega (X_1, \dots, X_{2k}) = s^k(X_1 \wedge \cdots \wedge
  X_{2k}) \omega_1 \cdots \omega_{2k+1}, \qquad \text {for} \quad X_i \in
  V_\un, \\    \Omega (X_1, \dots, X_{2k}) = 0, \qquad \text {if one} \quad X_i \in
  V_\ze.
\end{gather*}
If a $2k$-cocycle $\Omega'$ verifies $\Omega'(X_1, \dots, X_{2k}) =
\Omega(X_1, \dots, X_{2k}) $ for all $X_i \in V$, then $\Omega' =
\Omega$ ${\rm{ mod}} \ B^{2k}$.}

\end{enumerate}

As a consequence, $\Cc(2n{+}1, 2k)$ is rigid if $k {\neq} 1$. In the case
of $\Cc(2n{+}1,2)$, since $H^3 (\Cc(2n{+}1,2))$ $= \{ 0 \}$, there does
exist non-trivial deformations, and more precisely, a universal
deformation formula (see Appendix~\ref{appendix2}).

\medskip

\noindent
{\bf Theorem 5.}
{\it Let $\Ac_\Lambda$ be the $\CC[\Lambda]$-algebra generated by
  $V_\ze = \CC   P$ and $V_\un = \spa \{ E_+, E_- \}$ with relations:
  \begin{gather*}
   [E_+,E_-]_\bl = -\tfrac14 + \Lambda P, \qquad P^2 =1 \qquad \text{and}\qquad
  PE_\pm = -E_\pm P.
\end{gather*}

  Then $\Ac_\Lambda$ is a non-trivial polynomial deformation of
  $\Cc(1,2)$ and a universal deformation formula.  }

\medskip

This algebra $\Ac_\Lambda$ is a particular case (the simplest one) of
a symplectic ref\/lection algebra~\cite{Etingof-Ginzburg}.

Given an algebra $A$, any deformation $A_\Lambda$ of $A$ naturally
produces a deformation $\Mcc_r(A_\Lambda)$ of $\Mcc_r(A)$ and
conversely, any deformation of $\Mcc_r(A)$ is of type
$\Mcc_r(A_\Lambda)$ up to equivalence (see Appendix~\ref{appendix3}). Therefore,
Theorem~1 allows us to conclude that $\Mcc_{2^n}(\Ac_\Lambda)$ is a
universal deformation formula of $\Cc(2n+1,2)$.

Another presentation of this deformation is given in Section~\ref{section6}. We
introduce algebras $\agln$ and $\apln$ ($\Lambda$ formal, $\lambda \in
\CC, n \in \NN)$ by generators and relations, with $\agl(0) =
\agl$. We show that they can be constructed using Ore extensions, and
that they have a periodic behavior:

\medskip

\noindent
{\bf Periodicity Lemma 2.}
{\it $\agln \simeq \Cc(2n) \otimes \agl \simeq \Mcc_{2^n} (\agl)$ and the
same result holds for $\apl$.}

\medskip

In Sections \ref{section7} and~\ref{section8}, we describe $\agln$ and $\apln$ using the
enveloping algebra $\Uc$ of the Lie superalgebra $\osp(1,2)$ and its
primitive quotients \cite{P90}. Denoting by $\apl$ the algebra
$\apl(0)$, one has:

\medskip

\noindent
{\bf Theorem 6.}
\begin{enumerate}\itemsep=0pt

\item[$1.$] {\it $\agln \simeq \Mcc_{2^n} (\Sk_2 \ltimes \Uc)$ and $\apln \simeq
  \agln / (\Lambda-\lambda) \agln$.}

\item[$2.$] {\em Let $C$ be the Casimir element of $\Uc$ and $\Bc_c = \Uc /
    (C -c) \Uc$, $c \in \CC$. Then $\apl \simeq \Bc_{\lambda^2 -
      \frac{1}{16}}$ if $\lambda \neq 0$ and $\Ac_0 = \Cc(1,2)$.}

\item[$3.$] {\it If $\lambda^2 \neq \left(h + \frac14 \right)^2$, $2h \in \NN$,
  then $\apln$ is a simple algebra.

    If $\lambda^2 = \left(h + \frac14 \right)^2$, $2h \in
  \NN$, then $\apln$ is a primitive algebra with a unique non-zero
  two-sided ideal that is the kernel of an irreducible representation
  $\pi_h$ of dimension $2^{n} \left( 4h +1 \right)$.}

\item[$4.$] {\it $\agln$ is a FCR-algebra. Any finite-dimensional representation
  of $\apln$ is completely reducible and isotypical of type $\pi_h$,
  if $\lambda^2 = \left(h + \frac14 \right)^2$, $2h \in \NN$.}
\end{enumerate}

We obtain in this way all primitive quotients of $\agln$: $2^{2n}
\left( 4h +1 \right)^2$-dimensional quotients ($2 h \in \NN$) and
inf\/inite-dimensional ones: $\apln \simeq \Bc_{\lambda^2 -
  \frac{1}{16}}$ if $\lambda \neq 0$ and $\Ac_0(n) = \Cc(2n+1,2)$.

The isomorphism $\agln \simeq \Mcc_{2^n} (\Sk_2 \ltimes \Uc)$ is
useful to construct representations. Remark that representations of
$\Sk_2 \ltimes \Uc$ are merely graded representations of $\Uc$. Then,
from a graded representation of $\osp(1,2)$ on a space $H$, one
constructs a representation of $\agln$ on $H^{2^n}$. All irreducible
f\/inite-dimensional representations are obtained in this way: from
irreducible $(4h+1)$-dimensional $(2h \in \NN)$ representations of
$\osp(1,2)$, one obtains $\pi_h$. From the metaplectic (oscillator)
representation of $\osp(1,2)$, natural inf\/inite-dimensional
representations can be also obtained, using Dunkl-type formulas given
in~\cite{LP}.

Finally, we extend the obtained supersymmetry by showing that $\agl$
is a quotient of $\Sk_2 \ltimes \Uc (\osp(2,2))$. Therefore $\agln$ is
a quotient of $\Mcc_{2^n} \left( \Sk_2 \ltimes \Uc(\osp(2,2) \right)$.

There are three Appendices: the f\/irst one is a short introduction to
Hochschild cohomo\-logy relating it to deformation theory. In the second
Appendix, we explain, with proofs, what a~``universal deformation
formula'' is. We show in the third Appendix, that given an algebra
$A$, deformations of $\Mcc_r(A)$ are of type $\Mcc_r(A_\Lambda)$ up to
equivalence, where $A_\Lambda$ is a deformation of $A$. Results in
Appendices~\ref{appendix2} and~\ref{appendix3} are known, but since we have not found a convenient
reference, short proofs are given.

\section[Clifford algebras and Weyl algebras]{Clif\/ford algebras and Weyl algebras} \label{Section00}

We begin by recalling some classical properties of Clif\/ford and Weyl
algebras needed in the paper. Throughout this section, we denote
  by $[\cdot,\cdot]$ the super bracket and by
  $[\cdot,\cdot]_\bl$ the Lie bracket.

\subsection{Weyl algebras}

Let $k \in \NN^*$ and $\Sb_{2k} = \CC[p_1, q_1, \dots, p_k, q_k]$ be
the polynomial algebra in $2k$ indeterminates equipped with the
Poisson bracket:
\[ \{F,G\} = \sum_{i=1}^k \left( \dfrac{\partial F}{\partial p_i}
\dfrac{\partial G}{\partial q_i} - \dfrac{\partial F}{\partial q_i}
\dfrac{\partial G}{\partial p_i} \right), \qquad \forall \, F, G \in
\Sb_{2k}.
\]

Let $\wp \colon \Sb_{2k} \otimes \Sb_{2k} \to \Sb_{2k} \otimes
\Sb_{2k}$ be the operator def\/ined by:
\[\wp := \sum_{i=1}^k \left( \dfrac{\partial}{\partial p_i} \otimes
\dfrac{\partial}{\partial q_i} - \dfrac{\partial}{\partial q_i}
\otimes \dfrac{\partial}{\partial p_i} \right). \]

Let $m$ be the product of $\Sb_{2k}$ and $t$ be a formal parameter (or
$t \in \CC$). A new associative product $m_\tsti$ is def\/ined by:
\begin{equation}\label{1.2A}
m_\tsti := m \circ \expo\left( \frac{t}{2} \wp \right).
\end{equation}

This product $m_\tsti$ is a deformation of $\Sb_{2k}$ guided by the
Poisson bracket.

\begin{defn}
  The {\em Weyl algebra} $\Wb_{2k}$ is the vector space $\Sb_{2k}$
  endowed with the product $\sta:=m_\tstu$, called the {\em Moyal
    product}.
\end{defn}

A presentation of $\Wb_{2k}$ is given by generators $\{p_1,q_1, \dots,
p_k, q_k\}$ and relations:
\[ [u,v]_\bl = \{u,v \} \cdot 1, \qquad \forall \, u,v \in \spa\{p_1,q_1,
\dots, p_k, q_k\}. \]

Structurally, $\Wb_{2k}$ is central, simple, naturally $\ZZ_2$-graded
by the parity of $\Sb_{2k}$ and has a~supertrace~\cite{PU05}:
\begin{equation}\label{1.4A}
\Str(F):=F(0), \qquad \forall \, F \in \Wb_{2k}.
\end{equation}

The space $\Mb_k := \CC[x_1, \dots, x_k]$ is a faithful simple
$\Wb_{2k}$-module if we realize $p_i$ as $\frac{\partial}{\partial
  x_i}$ and~$q_i$ as the multiplication by~$x_i$, for all $i = 1,
\dots, k$. In the sequel, $\Mb_k$ is called the {\em metaplectic
  representation} of the Weyl algebra $\Wb_{2k}$.

The algebra of operators $\Lc(\Mb_k)$ appears as a completion of the
Weyl algebra: $\Wb_{2k}$ is the algebra of dif\/ferential operators of
f\/inite order, and any element $T$ in $\Lc(\Mb_{k})$ is a dif\/ferential
operator, in general of inf\/inite order (i.e., in the formal sense, the
sum is not f\/inite), given by the formula:
\begin{equation}\label{1.5A}
T = \sum_N \dfrac{1}{N!} \left(  m \circ (T \otimes \Sc) \circ \Delta
\big(x^N\big) \right)   \frac{\partial^N}{\partial x^N},
\end{equation}
where $\Sc$ is the antipode of $\Mb_{2k}$, $\Delta$ is its co-product,
$x^N := x_1^{n_1} \cdots x_k^{n_k}$ and $\frac{\partial^N}{\partial
  x^N} := \frac{\partial^{n_1 + \cdots + n_k}}{\partial x_1^{n_1} \cdots
  \partial x_k^{n_k}}$ if $N = (n_1, \dots, n_k)$. This formula gives
the (formal) symbol of $T$ in the normal ordering, and for
well-behaved $T$, its (formal) symbol in the Weyl ordering (see~\cite{PU05}).

From the point of view of deformation theory, $\Wb_{2k}$ is
rigid. More precisely, we have $H^r (\Wb_{2k})$ $= \{ 0 \}$, for all $r
>0$~\cite{Srid}.

We refer to \cite{PU05} for more details on the Weyl algebra in the
context of this section.

\subsection[Clifford algebras]{Clif\/ford algebras}

Let $n \in \NN^*$ and $\Wedge_n$ be the Grassmann algebra in $n$
anti-commutative variables $\omega_1, \dots, \omega_n$. Recall that
$\Wedge_n$ is $\ZZ$-graded. Denote by $\partial_1, \dots, \partial_n$
the super-derivations def\/ined by $\partial_i(\omega_j) = \delta_{ij}$,
$\forall \, i, j$. The algebra $\Wedge_n$ is endowed with a super
Poisson bracket:
\[\{\Omega, \Omega'\} = 2 (-1)^{\deg_\ZZ(\Omega)+1}
\sum_{i=1}^n \partial_i(\Omega) \wedge \partial_i(\Omega'), \]
for all
$\Omega, \Omega' \in \Wedge_n$ \cite{PU07}. We def\/ine the operator
$\wp$ of $\Wedge_n \otimes \Wedge_n$ by:
\[\wp := \sum_{i=1}^n \partial_i \otimes \partial_i,\]
where $\otimes$ is the graded tensor product of operators.

Let $m_\wedge$ be the product of $\Wedge_n$ and $t$ be a formal
parameter (or $t \in \CC$).  A new product $m_\tsti$ can be def\/ined by
(see \cite{PU07}):
\begin{equation}\label{2.2A}
m_\tsti := m_\wedge \circ \expo(-t \wp).
\end{equation}

\begin{defn}
The {\em Clifford algebra} $\Cc(n)$ is the vector space $\Wedge_n$
equipped with the product $\sta:= m_\tstu$.
\end{defn}

There is a $\ZZ_2$-gradation on $\Cc(n)$ def\/ined by $\dzd(\omega_i) =
1$, for all $i = 1, \dots, n$. A presentation of $\Cc(n)$ is given by
basic generators $\omega_1, \dots, \omega_n$ and relations:
\[[v,v']= \{v,v'\} \cdot 1, \qquad \forall \, v, v' \in \spa\{\omega_1,,
\dots, \omega_n \}.\]

In particular, we have:
\begin{gather*}
\omega_i^2 = 1, \quad \forall\, i, \qquad \omega_i \sta \omega_j +
\omega_j \sta \omega_i = 0, \quad \forall \, i \neq j \qquad \text{and}\\
\omega_{i_1} \wedge \cdots \wedge \omega_{i_p} = \omega_{i_1} \sta
\cdots \sta \omega_{i_p} \qquad \text{if} \quad i_1 < i_2 < \cdots < i_p, \quad p
\leq n.
\end{gather*}

\subsection[Even Clifford algebras]{Even Clif\/ford algebras} \label{eca}

For $i = 1, \dots, n$, let $P_i= \partial_i$ and $Q_i = x_i \wedge \
$. be respectively the operator of derivation and multiplication of
the Grassmann algebra $\Phi_n$ in $n$ anti-commutative variables $x_1,
\dots, x_n$. The operators $\omega_{2j-1} = Q_j + P_j$ and
$\omega_{2j} = i (Q_j - P_j)$, $j = 1, \dots, n$ verify the def\/ining
relations of the Clif\/ford algebra $\Cc(2n)$, so there is a
homomorphism from $\Cc(2n)$ onto the algebra of dif\/ferential operators
$\Diff(\Phi_n)$. It is easy to see that $\dim (\Cc(2n)) =
\dim(\Diff(\Phi_n))= \dim(\Lc(\Phi_n))$, so we can identify $\Cc(2n)=
\Diff(\Phi_n) = \Mcc_{2^n}(\CC)$, where $\Mcc_{2^n}(\CC)$ denotes the
algebra of complex matrices of order $2^n$. As a consequence, $\Phi_n$
is the unique simple $\Cc(2n)$-module, called the {\em spin
  representation} of $\Cc(2n)$.

Structurally, the even Clif\/ford algebra $\Cc(2n) = \Mcc_{2^n}(\CC)$ is
simple and its center is $\CC$.  From the point of view of deformation
theory, $\Cc(2n)$ is rigid and we have $H^r(\Cc(2n)) = \{0 \}$ if $r >
0$.

Since $\Cc(2n) = \Mcc_{2^n} (\CC)$, there is a natural trace on
$\Cc(2n)$ that can be written in an analogous way as in (\ref{1.4A}):
\[\Tr(\Omega):= 2^n   \Omega(0), \qquad \forall \, \Omega \in \Cc(2n).\]

There is also a similar formula to (\ref{1.5A}) in the case of
$\Cc(2n)$. In other words, any operator $T \in \Lc(\Phi_n)$ is
dif\/ferential and an explicit formula is given by:
\[T = \sum_{I \in \{0, 1\}^n} (-1)^{\theta(I,I)} \left( m_\wedge \circ
  (T \otimes \Sc) \circ \Delta\big(x^I\big) \right) \wedge \partial^I, \]
where $\theta$ is the bilinear form on $\NN^n$ associated to the
matrix $(\theta_{rs})_{r,s=1}^n$ with $\theta_{rs} = 1$ if $r>s$ and~$0$ otherwise, $\otimes$ is the non-graded tensor product of
operators, $\Sc$ is the antipode of $\Phi_n$, $\Delta$ is its
co-product, $x^I := x_1^{i_1} \wedge \cdots \wedge x_n^{i_n}$ and
$\partial^I := \partial_1^{i_1} \circ \cdots \circ \partial_n^{i_n}$ if
$I = (i_1, \dots, i_n)$.

\subsection[Periodicity of Clifford algebras]{Periodicity of Clif\/ford algebras} \label{1.4}

There is an algebra isomorphism between $\Cc(2n+k)$ and $\Cc(2n)
\otimes \Cc(k)$ since $\omega_1 \otimes 1, \dots, \omega_{2n} \otimes
1$ and $i^n  \omega_1 \sta \cdots \sta \omega_{2n} \otimes
\omega_j'$, $j =1, \dots, k$ in $\Cc(2n) \otimes \Cc(k)$ verify the
def\/ining relations of $\Cc(2n+k)$ thanks to the formula:
\[(\omega_1 \sta \cdots \sta \omega_{2n})^2 = (-1)^n.\]

It results that:
\[\Cc(2n) \simeq  \Cc(2)^{\otimes_n} \qquad \text{and} \qquad \Cc(2n+1) \simeq \Cc(2n) \otimes \Cc(1) \simeq \Mcc_{2^n} \left( \Cc(1) \right).\]

\subsection[Odd Clifford algebras]{Odd Clif\/ford algebras}

Since $\Cc(1) \simeq \CC \times \CC$, from the isomorphism $\Cc(2n+1)
\simeq \Mcc_{2^n} \left( \Cc(1) \right)$ it follows that $\Cc(2n+1)$
is the product $\Cc(2n) \times \Cc(2n)$. Therefore $H^r ( \Cc(2n+1)) =
\{0\}$ if $r>0$, and that implies that $\Cc(2n+1)$ is rigid.

We will make more explicit the above isomorphism: $\Cc(2n+1) \simeq
\Cc(2n) \times \Cc(2n)$.  The element $z=\omega_1 \sta \cdots \sta
\omega_{2n+1}$ is central and verif\/ies $z^2 = (-1)^n$. Set $Z =
\spa\{1,z\}$. Then $\Cc(2n+1) \simeq Z \otimes \Cc(2n)$ as
algebras. Let $z_+ = \frac12 (1 +i^nz)$ and $z_- = \frac12 (1 -
i^nz)$. Therefore $z_\pm^2 = z_\pm$, $z_+ \sta z_- = z_- \sta z_+ = 0$
and $1 = z_+ + z_-$. We conclude that
\[\Cc(2n+1) = z_+ \sta \Cc(2n) \oplus z_- \sta \Cc(2n),\]
that is, a reduction of $\Cc(2n+1)$ in a direct sum of two ideals
isomorphic to $\Cc(2n)$ as algebras. It follows that $\Cc(2n+1)$ has
exactly two $2^n$-dimensional simple modules built from the spin
representation of $\Cc(2n)$. To give more details, we need the
following lemma:

\begin{lem}
  Let $P$ be the natural parity of $\Phi_n$. Then in the spin
  representation of $\Cc(2n)$, one~has:
\begin{gather*}
\omega_1 \sta \cdots \sta \omega_{2n} = i^n   P.
\end{gather*}
\end{lem}

\begin{proof}
  We set $T =\omega_1 \sta \cdots \sta \omega_{2n}$. The operator $T$
  of $\Phi_n$ is diagonalizable since $T ^2 = (-1)^n$, and it commutes
  with $P$. We denote by $S_{\ze, {\pm i^n}}$ and
  $S_{\un, {\pm i^n}}$ its eigenspaces in $S_\ze$ and
  $S_\un$ respectively, where $S=\Phi_n$. The subspaces
  $S_{\ze, {i^n}} \oplus S_{\un, {-i^n}}$ and
  $S_{\ze, {-i^n}} \oplus S_{\un, {i^n}}$ are
  $\Cc(2n)$-stable since $T$ anti-commutes with $\omega_i$, $1 \leq i
  \leq 2n$. It follows that $T = \pm i^n P$. To determine the sign, we
  compute:
\begin{gather*}
T(1)  =  i^n (Q_1 + P_1) \sta (Q_2 - P_2)  \sta \cdots  \sta (Q_n + P_n)\sta(Q_n
- P_n)(1) = i^n.
\end{gather*}

Finally, we obtain $T = i^n P$.
\end{proof}

The element $z$ is central, $z^2 = (-1)^n$, hence $z= \pm i^n$ in any
simple $\Cc(2n+1)$-module. Since $z = (\omega_1 \sta \cdots \sta
\omega_{2n}) \sta \omega_{2n+1}$, using the lemma we obtain the spin
representations $\Phi^\pm_n$ of $\Cc(2n+1)$ as follows: $\Cc(2n)
\subset \Cc(2n+1)$ acts on $\Phi_n$ by the spin representation (see
Section~\ref{eca}), for $\Phi_n^+$, def\/ine $\omega_{2n+1}= P$ and for
$\Phi_n^-$, def\/ine $\omega_{2n+1}=-P$.

\subsection{Weyl algebras and supersymmetry}

Let $\Wb= \Wb_{2k}= \oplus_{r \geq 0} \Wb^{r}$ be the $\ZZ$-gradation
of the vector space $\Wb$. Recall that $\Wb$ is a $\ZZ_2$-graded
algebra. By (\ref{1.2A}), we have:
\[ [F,G]_\bl = \{F,G\}, \qquad \forall\, F \in \bigoplus_{r \leq 2} \Wb^r.\]

Let $\hk = \hk_\un \oplus \hk_\ze$ where $\hk_\un = \Wb^1= \spa\{p_i,
q_i, \, i = 1, \dots, k\}$ and $\hk_\ze = \Wb^2$. The super bracket
stabilizes $\hk$. Besides, $\hk$ is isomorphic to the Lie superalgebra
$\osp(1,2k)$. In particular, $\hk_\ze \simeq \spk(2k)$ and the adjoint
action of $\hk_\ze$ on $\hk_\un$ is the standard action of $\spk(2k)$
on $\CC^{2k}$. As a consequence, $\Wb$ is a semisimple
$\hk_\ze$-module for the adjoint action and $\Wb= \oplus_{r \geq 0}
\Wb^{r}$ is its reduction in isotypical components.

By (\ref{1.2A}), we have:
\[[v,F] = 2 v F, \quad \forall \, F \in \Wb^{2r+1} \qquad \text{and} \qquad [v,F]_\bl
= \{v, F\}, \quad \forall \,F \in \Wb^{2r}, \quad v \in \Wb^1.
\]

Therefore $\Wb$ is also semi-simple for the adjoint action of $\hk$
and $\Wb= \oplus_{r \geq 0} \Ab^{r}$ is its reduction into isotypical
components, where $A_0 =\CC$ and $\Ab^r = \Wb^{2r-1} \oplus \Wb^{2r}$,
if $r>0$. We refer to~\cite{Musson} or~\cite{PU05} for more details.

\subsection[Clifford algebras and symmetry]{Clif\/ford algebras and symmetry}

Let $\Cc = \Cc(n)$. There is a $\ZZ$-gradation on the vector space
$\Cc$ and, as an algebra, $\Cc$ is $\ZZ_2$-graded. By (\ref{2.2A}), we
have:
\[ [\Omega, \Omega'] = \{\Omega,\Omega'\}, \qquad \forall\, \Omega \in
\bigoplus_{r \leq 2} \Cc^r.\]

Let $\g = \g_\un \oplus \g_\ze$ where $\g_\un = \Cc^1$ and $\g_\ze =
\Cc^2$. The Lie bracket stabilizes $\g$. Moreover, $\g$ is isomorphic
to the Lie algebra $\ok(n+1)$. In particular, $\g_\ze \simeq \ok(n)$
and the adjoint action of $\g_\ze$ on~$\g_\un$ is the standard action
of $\ok(n)$ on $\CC^{n}$. The direct sum $\g = \g_\un \oplus \g_\ze$
is a $\ZZ_2$-gradation for the Lie algebra $\g$, that is $[\g_i,
\g_j]_\bl \subset \g_{i+j}$ (this is not a graded Lie algebra!). For
the adjoint action, $\Cc$ is a semisimple $\g_\ze$-module and is
isomorphic to the $\ok(n)$-module $\Wedge \CC^n$, whose reduction into
isotypical components is well-known (see~\cite{Samelson} or~\cite{Fulton-Harris}). The reduction into isotypical components of the
$\g$-module $\Cc$ can be deduced, but this is simply not the subject
of this paper.

\section[Clifford-Weyl algebras]{Clif\/ford--Weyl algebras} \label{Section01}

We recall the construction of the exterior algebra of a $\ZZ_2$-graded
vector space $V = V_\ze \oplus V_\un$: let $\Wedge := \Wedge V_\ze$ be
the exterior algebra of $V_\ze$ and $\Sb:= \sym(V_\un)$ be the
symmetric algebra of $V_\un$. Using their $\ZZ$-gradation, def\/ine a
$\ZZ \times \ZZ_2$ gradation on $\Wedge$ and on $\Sb$ by
\[\Wedge{}^{(i,\zero)} = \Wedge{}^i, \quad \Wedge{}^{(i,\um)} = \{0\} \qquad \text{and}
\qquad \Sb^{(i,\overline{i})} =
\Sb^i,\quad  \Sb^{(i,\overline{j})} = \{0\} \quad \text{if} \quad \overline{i}
\neq \overline{j}.\]

The exterior algebra of $V$ is the $\ZZ \times \ZZ_2$-graded algebra
\[\Es := \Wedge \underset{\ZZ \times \ZZ_2}{\otimes} \Sb= \Wedge
\underset{\ZZ}{\otimes} \ \Sb\]
endowed with the product:
\[ (\Omega \otimes F) \wedge (\Omega' \otimes F') = (-1)^{f\omega'} (\Omega \wedge \Omega') \otimes F F', \]
for all $\Omega \in \Wedge$, $\Omega' \in \Wedge^{\omega'}$, $F \in
\Sb^f$, $F' \in \Sb$. We have
\[A' \wedge A = (-1)^{aa' + bb'} A \wedge A', \qquad  \forall \, A \in
\Es^{(a, \overline{b})}, \quad A' \in \Es^{(a', \overline{b'})}\]
and that means that $\Es$ is $\ZZ \times \ZZ_2$-commutative.

Now, assume that the dimension of $V_\un$ is even, say $\dim(V_\un ) =
2k$. Set $n = \dim( V_\ze)$.

We have def\/ined Poisson brackets on $\Wedge$ and $\Sb$ in Section~\ref{Section00}. Associated operators $\wp_\lw$ and $\wp_\lsb$ are
respectively def\/ined on $\Wedge \otimes \Wedge$ and $\Sb \otimes \Sb$.

A super $\ZZ \times \ZZ_2$-Poisson bracket on $\Es$ is then def\/ined
by:
\[ \{ \Omega \otimes F, \Omega' \otimes F' \} = (-1)^{f \omega'}
\left( \{ \Omega, \Omega'\} \otimes F F' + (\Omega \wedge \Omega')
  \otimes \{F,F' \} \right),\] for all $\Omega \in \Wedge$, $\Omega'
\in \Wedge^{\omega'}$, $F \in \Sb^f$, $F' \in \Sb$. Now, let
$\sigma_{23}$ and $\wp$ be operators on $\Es \otimes \Es$ def\/ined by:
\begin{gather*}
 \sigma_{23} (\Omega \otimes F \otimes \Omega' \otimes F') = (-1)^{f
    \omega'} \Omega \otimes \Omega' \otimes F \otimes F',
    \\
    \wp =
  \sigma_{23} \circ (-2 \wp_\lw \otimes \Id + \Id \otimes \wp_\lsb)
  \circ \sigma_{23},
\end{gather*}
for all $\Omega \in \Wedge$, $\Omega' \in \Wedge^{\omega'}$, $F \in
\Sb^f$, $F' \in \Sb$.

Let $t$ be a formal parameter (or $t \in \CC$). A new product
$m_\tsti$ on $\Es$ is def\/ined from these operators and from the
product $m_{{\Es}}$ on $\Es$ by:
\begin{equation} \label{pst}
m_\tsti := m_{{\Es}} \circ \expo \left(\frac{t}{2} \wp \right).
\end{equation}

Since $m_\tsti = m_\lw \otimes m_\lsb\circ \expo \left( -t \wp_\lw
\right) \otimes \expo \left( \frac{t}{2} \wp_\lsb \right)\circ
\sigma_{23}$, it results that $m_\tsti$ is exactly the $\ZZ_2 \times
\ZZ_2$-graded tensor algebra product:
\[ \Cc^t(n) \otimes_{\ZZ_2 \times \ZZ_2} \Wb_{2k}^t = \Cc^t(n)
\otimes_{\ZZ_2} \Wb_{2k}^t,
\] where $ \otimes_{\ZZ_2}$ means the graded
tensor product with respect to left $\ZZ_2$-gradations, $\Cc(n)^t$
denotes the algebra equipped with product $m_\tsti$ and similarly for
$\Wb_{2k}^t$ (see Section~\ref{Section00}). By def\/inition, $m_\tsti$
is a deformation of $ m_{{\Es}}$ guided by the Poisson
super bracket.

\begin{defn}
The {\em Clifford--Weyl algebra} $\Cc(n,2k)$ is the vector space $\Es$
endowed with the product $\sta := m_\tstu$.
\end{defn}

Denote by $\{ \omega_1, \dots, \omega_n\}$ and $\{p_1,q_1, \dots,
p_k,q_k\}$ respectively the basis of $\Cc(n)$ and $\Wb_{2k}$ as in
Section~\ref{Section00}. The algebra $\Cc(n,2k)$ has a presentation
given by generators $\{\omega_1, \dots, \omega_n, p_1,q_1, \dots,
$ $p_k,q_k\}$ and relations
\begin{gather*}
    [\omega_i, \omega_j]_+ = 2 \delta_{ij}, \quad [p_i, q_j]_- =
  \delta_{ij}, \quad [p_i, p_j]_- = [q_i,q_j]_- =0 \quad \text{if} \quad i \neq
  j \\   \text{and} \quad  [\omega_i, p_j]_+ = 0, \quad [\omega_i, q_j]_+ = 0, \quad \forall \, i, j,
\end{gather*}
where $[A,B]_\pm := A \sta B \pm B \sta A$.

When $n = 2 \ell$ is even, we set
\[P_j = \tfrac12 ( \omega_{2j-1}+i \omega_{2j}) \qquad \text{and} \qquad Q_j =
\tfrac12 ( \omega_{2j-1}-i \omega_{2j}) \qquad \text{for} \quad j \leq \ell.\]

The f\/irst two relations above become:
\[ [P_i,Q_j]_+ = \delta_{ij}, \qquad  [P_i, P_j]_+ = [Q_i,Q_j]_+ =0.\]

Consider now $\Phi_\ell$ the Grassmann algebra in $\ell$
anti-commutative variables $\xi_1, \dots, \xi_\ell$, $\Mb_k$ the
polynomial algebra in $k$ variables $x_1, \dots, x_k$ and
$\SM(\ell,k)$ the exterior algebra of the $\ZZ_2$-graded space $W =
W_\ze \oplus W_\un$ with $W_\ze = \spa \{ \xi_1, \dots, \xi_\ell \}$
and $W_\un = \spa \{x_1, \dots, x_k \}$. There is a~$\Cc(2 \ell, 2
k)$-module structure on $\SM(\ell, k)$ given by $P_i =
\frac{\partial}{\partial \xi_i}$, $ Q_i = \xi_i \wedge \cdot \, $
($1 \leq i \leq \ell$), $p_j = \frac{\partial}{\partial x_j}$, $ q_j =
x_j \ \cdot \ $ ($1 \leq j \leq k$). Besides, $\SM(\ell, k)$ is a
simple $\Cc(2 \ell, 2 k)$-module. In the sequel, we call $\SM(\ell,
k)$ the {\em spin-metaplectic representation} of $\Cc(n,2k)$. This
provides a homomorphism from $\Cc(2 \ell,2k)$ onto the algebra of
dif\/ferential operators of the $\ZZ_2 \times \ZZ_2$-graded exterior
algebra $\SM(\ell,k)$. We will show later that $\Cc(n,2k)$ is simple,
so we have actually an isomorphism, $\Cc(2 \ell, 2k) \simeq
\Diff(\SM(\ell,k))$ and that generalizes the cases of $\Cc(2 \ell)$
and $\Wb_{2k}$ seen in Section \ref{Section00}.

When $n = 2 \ell +1$ is odd, we obtain two spin-metaplectic
representations $\SM(\ell,k)^\pm$ of $\Cc(2 \ell +1$, $2k)$ by acting
$\Cc(2 \ell, 2k)$ on $\SM(\ell,k)$ as above and by setting $\omega_{2
  \ell +1} = Q$ for $ \SM(\ell,k)^+$ and $\omega_{2\ell +1} = -Q$ for
$\SM(\ell,k)^-$, where $Q$ is the parity:
\[ Q(\omega \otimes f) = (-1)^{\deg_\ZZ(\omega) +\deg_\ZZ(f)}  \omega \otimes f, \qquad \forall \, \omega \in \Phi_\ell, \quad f \in \Mb_k.\]

It will be shown later that if $k \neq 0$, $\Cc(2 \ell +1, 2 k)$ is
simple and as a consequence, both representations $\SM(\ell, k)^\pm$
are faithful.

\section[Periodicity of Clifford-Weyl algebras]{Periodicity of Clif\/ford--Weyl algebras}\label{section3}

Clif\/ford algebras have a periodic behavior Section~\ref{1.4} and we now show
that this periodicity can be extended to Clif\/ford--Weyl algebras. We
denote by $\Cc(r)$, the Clif\/ford algebra in $r$ variables and by
$\Wb_{2k}$, the Weyl algebra constructed from $2k$ variables (see
Section \ref{Section00}).

\begin{lem}[Periodicity Lemma 1] \label{3.1}
\[\Cc(2m+n, 2k)   \simeq \Cc(2m) \otimes \Cc(n,2k).
\]
\end{lem}

\begin{proof}
  Let $\{ \omega_1, \dots, \omega_{2m}\}$ and $\{\omega_1', \dots,
  \omega_{n}',p_1,q_1, \dots, p_k,q_k\}$ be respectively the set of
  generators of $\Cc(2m)$ and $\Cc(n,2k)$. Let $z = i^m \omega_1 \sta
  \cdots \sta \omega_{2m}$. So $z^2 =1$ and $z$ anti-commutes with
  $\omega_1, \dots, \omega_{2m}$. The following elements of $\Cc(2m)
  \otimes \Cc(n,2k)$: $\omega_1 \otimes 1, \dots, \omega_{2m} \otimes
  1$, $z \otimes \omega_1', \dots, z \otimes \omega_{n}'$, $z \otimes
  p_1, \dots, z \otimes p_k$ and $z \otimes q_1, \dots, z \otimes q_k$
  verify the def\/ining relations of $\Cc(2m+n,2k)$. Since they generate
  $\Cc(2m) \otimes \Cc(2n,k)$ as an algebra, we get an algebra
  homomorphism from $\Cc(2m+n,2k)$ onto $\Cc(2m) \otimes \Cc(n,2k)$.

  Denote by $\{\widetilde{\omega_1}, \dots, \widetilde{\omega}_{2m+n},
  \widetilde{p_1},\widetilde{q_1}, \dots,
  \widetilde{p_k},\widetilde{q_k}\}$ the set of generators of
  $\Cc(2m+n,2k)$. Let $\widetilde{z} =  i^m
  \widetilde{\omega_1} \sta \cdots \sta \widetilde{\omega}_{2m}$. So
  $\widetilde{z}^2 =1$, $\widetilde{z}$ anti-commutes with
  $\widetilde{\omega_1}, \dots, \widetilde{\omega}_{2m+n}$ and
  $\widetilde{z}$ commutes with $p_1, \dots, p_k$ and $q_1, \dots,
  q_k$. The following elements of $\Cc(2m+n,2k)$:
  $\widetilde{\omega_1}, \dots, \widetilde{\omega}_{2m}, \widetilde{z}
  \sta \widetilde{\omega}_{2m+1}, \dots$, $\widetilde{z} \sta
  \widetilde{\omega}_{2m+n}, \widetilde{z} \sta \widetilde{p_1},$
  $\widetilde{z} \sta \widetilde{q_1}, \dots, \widetilde{z} \sta
  \widetilde{p_k},\widetilde{z} \sta \widetilde{q_k}$ verify the
  def\/ining relations of $\Cc(2m) \otimes \Cc(n,2k)$, so we get the
  inverse homomorphism.
\end{proof}

\begin{cor} \label{isoscw} One has:
\begin{gather*}
\Cc(2n,2k) \simeq \Cc(2n) \otimes \Wb_{2k} \simeq \Mcc_{2^n} \left(
\Wb_{2k} \right) \qquad \text{and} \\
\Cc(2n+1,2k) \simeq \Cc(2n) \otimes \Cc(1,2k) \simeq \Mcc_{2^n} \left(
  \Cc(1,2k) \right),
\end{gather*}
 where $\Mcc_r(A)$ denotes the $r \times
r$-matrix algebra with coefficients in an algebra $A$ for a given $r
\in \NN^*$.
\end{cor}

{\samepage

\begin{cor} \label{2.5} \quad
\begin{enumerate}\itemsep=0pt

\item[$1.$] $\Cc(2n,2k)$ is simple with center $\CC$.

\item[$2.$] If $k \neq 0$, then $\Cc(2n+1,2k)$ is simple with center $\CC$.

\end{enumerate}
\end{cor}}

\begin{proof}
1. $\Cc(2n,2k) \simeq \Mcc_{2^n}(\Wb_{2k})$ is simple since
  $\Wb_{2k}$ is simple.

2. Since $\Cc(2n+1,2k) \simeq \Mcc_{2^n}(\Cc(1,2k))$, it is enough
  to prove the result for $\Cc(1,2k)$.

  But $\Cc(1,2k) \simeq \Sk_2 \ltimes \Wb_{2k}$ and $\Wb_{2k}$ is
  simple, so the result is a particular case of a general theorem in
  \cite{Montgomery}.

  For the sake of completeness, here is a direct proof: we write
  $\Cc(1,2k) = \Cc(1) \zdtimes \Wb_{2k}$ where $\Cc(1)$ is the
  Clif\/ford algebra generated by $\Ps$ such that $\Ps^2 = 1$. Recall
  that using the Moyal $\sta$-pro\-duct, the Weyl algebra $\Wb=
  \Wb_{2k}$ can be realized as a deformation of the polynomial algebra
  $\CC[p_1,q_1, \dots , p_k, q_k]$. Fix $p = p_1$ and $q = q_1$.

We have $[p, f]_\bl= \frac{\partial f}{\partial q}$, $\forall\, f \in
\Wb$.  In addition, for all $g \in \Wb$:
\begin{gather*}
  [p, \Ps \sta \ g]_\bl  =  p \sta \Ps \sta \ g - \Ps \sta \ g \sta p =
  - \Ps \sta (p\sta g + g \sta p)\\
\phantom{[p, \Ps \sta \ g]_\bl}{} = - \Ps \sta \left( p g +
    \tfrac12\{p,g\} + g p + \tfrac12 \{g,p \} \right) = -2 \Ps \sta (p g).
\end{gather*}

Let $I$ be a non-zero two-sided ideal of $\Cc(1,2k)$ and let $f+ \Ps
\sta g \in I$, $f+ \Ps \sta g \neq 0$. Then $[p, f+\Ps \sta g]_\bl\in
I$ gives $\frac{\partial f}{\partial q} - 2\Ps \sta (p g) \in I$ and
we can reiterate. Hence:
\begin{itemize}\itemsep=0pt
\item if $g =0$, then $f \in I$. It follows that $I \cap \Wb \neq
  \{0\}$;

\item if $g \neq 0$, since there exists $j$ such that
  $\frac{\partial^j f}{\partial q^j} = 0$, one has $(-1)^j 2^j \Ps
  \sta (p^j g) \in I$, implying $p^j g \in I$. But $p^j g \neq 0$, so
  it follows that $I \cap \Wb \neq \{0\}$ as well.
\end{itemize}

In both cases, $I \cap \Wb$ is a non-zero ideal of the Weyl algebra
$\Wb$. Since $\Wb$ is simple, $I \cap \Wb = \Wb$. So $1 \in I$ and we
conclude that $I= \Cc(1,2k)$.

The center of $\Cc(1,2k)$ is $\CC$ since the center of $\Wb$ is $\CC$.
\end{proof}

\begin{rem} \label{3.4}
  Let us f\/irst recall what a smash product is. Let $A$ be an algebra,
  $G$ a f\/inite group acting on $A$ by automorphisms and $\CC[G]$ the
  group algebra. The {\em smash product} $G \ltimes A$ is the algebra
  with underlying space $A \otimes \CC[G]$ and product def\/ined by:
  \[
  (a \otimes g) (a' \otimes g') = a(g.a') \otimes gg', \qquad \forall \,
  a,a' \in A, \quad g,g' \in G.
  \] Smash products used in this paper are
  def\/ined from a $\ZZ_2$-graded algebra $A$ and $G = \Sk_2 =
  \{-1,1\}$. Denoting by $\Ps$ the parity operator of $A$,
  $\Ps(a):=(-1)^{\deg(a)} a$, $\forall\, a \in A$, $\Sk_2$ acts on $A$
  by $(-1)^\alpha \cdot a := \Ps^\alpha(a)$, $\forall\, a \in A$, $\alpha =
  0,1$ and there is a corresponding smash product $\Sk_2 \ltimes
  A$. It is the algebra generated by $\Ps$ and $A$ with relations $\Ps
  a = \Ps(a)\Ps$, $\forall\, a \in A$ and $\Ps^2 = 1$. It is easy to
  check that $\Sk_2 \ltimes A$-modules and $\ZZ_2$-graded $A$-modules
  are exactly the same notion.

  Now consider the Clif\/ford--Weyl algebra $\Cc(1,2k)$. Using the
  $\ZZ_2$-graded structure of $\Wb_{2k}$, $\Cc(1,2k) \simeq \Sk_{2}
  \ltimes \Wb_{2k}$. Also $\Cc(1,2k)$ is isomorphic to a subalgebra of
  $\Mcc_{2}(\Wb_{2k})$:
  \[ \Cc(1,2k) \simeq \left\{ \begin{pmatrix} a & b \\ \sigma(b) &
      \sigma(a) \end{pmatrix}, \, a, b \in \Wb_{2k} \right\},
       \] where
  $\sigma$ is the parity operator of $\Wb_{2k}$.  In this isomorphism,
  $\omega_1 \in \Cc(1)$ is realized as the matrix $\begin{pmatrix} 0 &
    1 \\ 1 & 0 \end{pmatrix}$ and $\Wb_{2k}$ as
  $\left\{\begin{pmatrix} a & 0 \\ 0 & \sigma(a) \end{pmatrix}, \, a \in
    \Wb_{2k} \right\}$.

  Finally, $\Cc(1,2k)$ is isomorphic to the algebra generated by the
  parity operator $\Ps$ of $\Mb_k= \CC[x_1, \dots, x_k]$ and
  $\Wb_{2k}$, realized as the algebra of dif\/ferential operators of
  $\Mb_k$ (see Section~\ref{Section00}).
\end{rem}

\section[Clifford-Weyl algebras and supersymmetry]{Clif\/ford--Weyl algebras and supersymmetry}\label{section4}

Let us consider the $\ZZ_2 \times \ZZ_2$-graded algebra $\Cc(n,2k)$
and the subspace $V = V_\ze \oplus V_\un$ where $V_\ze
=\Cc(n,2k)_{(\un,\ze)}  = \Wedge^1_n$ and $V_\un = \Cc(n,2k)_{(\un,
  \un)} = \Sb^1_{2k}$ (see Section~\ref{Section01} for the notation).

If $k = 0$, then $V_\ze \oplus [V_\ze, V_\ze]_\bl$ is a Lie algebra
for the natural Lie bracket of the Clif\/ford algebra, isomorphic to
$\ok(n+1)$ and $[V_\ze, V_\ze]_\bl$ is a Lie subalgebra isomorphic to
$\ok(n)$ (for details, see Section~\ref{Section00}).

If $n=0$, then $V_\un \oplus [V_\un, V_\un]$ is a Lie superalgebra for
the natural super bracket of the Weyl algebra, isomorphic to
$\osp(1,2k)$ and $[V_\un, V_\un]$ is a Lie algebra isomorphic to
$\spk(2k)$.

To generalize this situation, we need some notation: for an element $a
\in \Cc(n,2k)$, denote its $\ZZ_2\times \ZZ_2$-degree by $\Delta(a) :=
(\Delta_1(a), \Delta_2(a))$. {\em We consider $\Cc(n,2k)$ as an
  algebra $\ZZ_2$-graded by $\Delta_2$ and we denote by $[\cdot,
    \cdot]$ the associated super bracket.}

  The proposition below shows how to realize $\osp(n+1,2k)$ as a Lie
   sub-superalgebra of $\Cc(n,2k)$. This important result
  was f\/irst obtained by \cite{Palev} for $\osp(2 \ell + 1, 2k)$. We
  propose here another method to show the same result, inspired by
  \cite{PU05} and based on a well-chosen twisted adjoint action.

\begin{prop}[\cite{Palev}] \label{4.1}
Let $\g = V \oplus [V,V]$. Then $\g$ is a Lie sub-superalgebra of
$\Cc(n,2k)$ isomorphic to $\osp(n+1,2k)$. Moreover
\[\gO =V_\ze \oplus [V_\ze, V_\ze] \oplus [V_\un, V_\un]\]
with $[V_\ze, V_\ze] \simeq \ok(n)$, $[V_\un, V_\un] \simeq \spk(2k)$,
$V_\ze \oplus [V_\ze, V_\ze] \simeq \ok(n+1)$ and $\gO \simeq
\ok(n+1)\times \spk(2k)$. Also,
\[\gI = V_\un \oplus [V_\ze, V_\un]\]
and $V_\un \oplus [V_\un, V_\un] \simeq \osp(1,2k)$. If we set $\hk =
[V_\ze, V_\ze] \oplus [V_\un, V_\un] \oplus [V_\ze, V_\un]$, then
$\hk \simeq \osp(n,2k)$.
\end{prop}

\begin{proof}
  By a case by case straightforward computation, using the product
  formula (\ref{pst}), we get the formula:
\begin{equation} \label{seis}
[[X,Y],Z] = 2 \left( \{Y,Z \} X - (-1)^{\Delta_2(X)\Delta_2(Y)} \{X,Z \} Y \right), \qquad \forall \,  X,Y,Z \in V,
\end{equation}
where $\{\cdot, \cdot\}$ is the super Poisson bracket
def\/ined in Section~\ref{Section01}.

Hence $[[V,V],V] \subset V$. If $H \in [V,V]$ and $X, Y \in V$, then:
\[[H,[X,Y]] = [[H,X],Y] + (-1)^{\Delta_2(H) \Delta_2(X)} [X, [H,Y]].\]

Using (\ref{seis}), we conclude that $[[V,V],[V,V]] \subset [V,V]$,
therefore $\g$ is a Lie superalgebra and $\hk$ is a sub-superalgebra.

To prove the isomorphisms, we set $V' = \CC \oplus V$. Def\/ine a
non-degenerate supersymmetric 2-form $(\cdot | \cdot)$ on
$V'$ by:
\[ (X|Y) := \{X,Y\}, \quad \forall \, X, Y \in V \qquad \text{and} \qquad (1|1) = -2.\]

Then formula (\ref{seis}) becomes parastatistics relations:
\begin{equation} \label{six}\tag{PS} [[X,Y],Z] = 2 \left( (Y|Z) X -
    (-1)^{\Delta_2(X)\Delta_2(Y)} (X|Z) Y \right), \qquad \forall \, X, Y, Z \in V.
\end{equation}

Next, we def\/ine the $\Delta_1$-twisted adjoint representation of the
Lie superalgebra $\Cc(n,2k)$:
\[ \ad'(a)(b) := a \sta b - (-1)^{\Delta_2(a)\Delta_2(b)+ \Delta_1(a)} b \sta a, \qquad \forall \, a, b \in \Cc(n,2k).\]

It is easy to check that it is indeed a representation. If $H \in
\hk$, $\ad'(H) = \ad(H)$, writing $H = [X,Y]$ and using (\ref{six}),
one obtains:
\[ (\ad'(H)(Z)|T) = -(-1)^{\Delta_2(Z) \Delta_2(H)}(Z| \ad'(H)(T)), \qquad
\forall \, T \in V,\]
henceforth $\ad'(\hk)(V') \subset V'$ and $\hk \subset
\osp(n,2k)$. Since both spaces have the same dimension $\frac{n
  (n-1)}{2} + 2nk+k(2k+1)$ (see \cite{Sch}), it follows $\hk \simeq
\osp(n,2k)$.

It remains to examine the action of $\ad'(X)$ on $V'$ when $X \in
V$. We have  $\ad'(X)(Y) = 0$ if $X \in V_\ib$, $Y \in
V_\jb$ with $\overline{i} \neq \overline{j}$. Moreover, if $X, Y \in
V_\ze$, then $\ad'(X)(Y) = X \sta Y + Y \sta X = \{ X, Y \} \cdot 1 =
(X|Y)$. If $X, Y \in V_\un$, then $\ad'(X)(Y) = X \sta Y - Y \sta X =
\{ X, Y \} \cdot 1 = (X|Y)$. Since $\ad'(X)(1) = 2X$, f\/inally
$(\ad'(X)(Y)|1) =$  $-2 (X|Y) $$=
-(-1)^{\Delta_2(X)\Delta_2(Y)} (Y|\ad'(X)(1))$. So $\g \subset
\osp(n+1,2k)$ and both spaces have the same dimension.
\end{proof}

\begin{cor}\label{4.2}
  Let $V = V_\ze \oplus V_\un$ be a $\ZZ_2$-graded space with
  $\dim(V_\ze)= n$ and $\dim(V_\un) = 2k$. Assume that $V$ is equipped
  with a non-degenerate supersymmetric bilinear form $(\cdot|
  \cdot)$. Let $A$ be the $\ZZ_2$-graded algebra generated by $V
  = V_\ze \oplus V_\un$ and relations \eqref{six}. Then $A$ is
  isomorphic to the enveloping algebra $\Uc(\osp(n+1,2k))$.
\end{cor}

\begin{proof}
  We denote by $[\cdot, \cdot]_A$ the super bracket of
  $A$. Proceeding exactly as in the proof of Proposition~\ref{4.1}, we
  show that $V + [V,V]_A$ is a Lie superalgebra using the
  parastatistics relations~(\ref{six}). From the def\/inition of $A$
  together with Proposition~\ref{4.1}, there is an algebra
  homomorphism from~$A$ onto~$\Cc(n,2k)$ that is the identity when
  restricted to $V$. This homomorphism induces a Lie superalgebra
  homomorphism from $V + [V,V]_A$ onto $V \oplus [V,V]$ (realized in
  $\Cc(n,2k)$ and isomorphic to $\osp(n+1,2k)$ by Proposition~\ref{4.1}). That implies $\dim(V + [V,V]_A) \geq
  \dim(\osp(n+1,2k))$.

  On the other hand, $\dim(V + [V,V]_A) \leq \dim(V \oplus [V,V])$
  since we can write
\[[V,V]_A = [V_\ze,V_\ze]_A + [V_\un,V_\un]_A +   [V_\ze,V_\un]_A\]
and $\dim([V,V]) = \dim(V_\ze \wedge V_\ze) + \dim(V_\ze \otimes
V_\un) + \dim(V_\ze V_\un)$.

It results that $V \oplus [V,V]_A \simeq \osp(n+1,2k)$. Remark that
the parastatistics relations hold in the enveloping algebra
$\Uc(\osp(n+1,2k))$ since they hold in $\osp(n+1,2k)$. To f\/inish, we
apply the universal property of $\Uc(\osp(n+1,2k))$.
\end{proof}

\begin{rem}
  The result in Proposition \ref{4.1} is helpful to obtain explicit
  descriptions of $\osp(n+1$, $2k)$ (for instance, the root system).
\end{rem}

\begin{rem}
  As observed in \cite{Palev}, the fact that generators of $\Cc(n)$
  (Fermi-type operators) and those of $\Wb_{2k}$ (Bose-type operators)
  anti-commute in $\Cc(n,2k)$ is a main argument to prove that the Lie
  sub-superalgebra that they generate is $\osp(n+1,2k)$. However, the
  periodicity of Clif\/ford--Weyl algebras, namely $\Cc(2n, 2k) \simeq
  \Cc(2n) \otimes \Wb_{2k}$, shows that it is always possible to
  obtain $\Cc(2n,2k)$ from commuting Bose-type and Fermi-type
  operators (that will not live in the Lie superalgebra
  $\osp(2n+1,2k)$, but rather in $\ok(2n+1) \times \osp(1,2k)$).
\end{rem}

{\em In the sequel, all $\star$ products will simply be denoted by
  juxtaposition.}

\section[Cohomology of Clifford-Weyl algebras]{Cohomology of Clif\/ford--Weyl algebras}\label{section5}

In Appendix~\ref{appendix1}, the reader can f\/ind a short introduction to Hochschild
cohomology of an algebra with coef\/f\/icients in itself.

By Periodicity Lemma~1 and Corollary \ref{isoscw}, we have
\[\Cc(2n,2k) \simeq \Mcc_{2^n} \left( \Wb_{2k} \right) \qquad \text{and} \qquad \Cc(2n+1,2k) \simeq \Mcc_{2^n} \left( \Cc(1,2k) \right).\]

But for an algebra $A$, $\Mcc_r(A)$ and $A$ have isomorphic cohomology
spaces \cite{Loday}, so it results that the cohomology of
Clif\/ford--Weyl algebras can be computed from the cohomology of
$\Wb_{2k}$ and $\Cc(1,2k)$:

{\samepage

\begin{prop}\quad
\begin{enumerate}\itemsep=0pt
\item[$1.$] $H^\ell(\Cc(2n,2k)) =  \{0\}$ if $\ell  > 0$.

\item[$2.$] $H^\ell(\Cc(2n+1,2k)) = H^\ell(\Cc(1,2k))$, for all $\ell$.
\end{enumerate}
\end{prop}}

\begin{proof}
It is enough to remark that $H^\ell(\Wb_{2k}) = \{ 0 \}$ if $\ell >
0$ \cite{Srid}.
\end{proof}

We now give more details on the identif\/ications in the above
Proposition. We use the isomorphisms in Corollary \ref{isoscw}:
$\Cc(2n,2k) \simeq \Cc(2n) \otimes \Wb_{2k}$ and $\Cc(2n+1,2k) \simeq
\Cc(2n) \otimes \Cc(1,2k)$. The letter $A$ denotes either $\Wb_{2k}$
or $\Cc(1,2k)$.

Since $\Cc(2n)$ is separable, we compute the cohomology of $\Cc(2n)
\otimes A$ using normalized $\Cc(2n)$-relative cochains (see
\cite{GS}), that is, cochains
\[\Omega : (\Cc(2n) \otimes A)^\ell \to \Cc(2n) \otimes A\]
that verify:
\begin{gather*}
   \Omega(C  a_1, a_2, \dots, a_\ell) = C  \Omega(a_1, \dots,
  a_\ell) ,
  \\
   \Omega(a_1,\dots,   a_iC, a_{i+1}, \dots, a_\ell)
  = \Omega (a_1, \dots, a_i, C  a_{i+1}, \dots,
  a_\ell), 
  \\  \Omega(a_1, \dots, a_\ell  C) = \Omega(a_1,
  \dots, a_\ell)  C , \\
  \Omega(a_1, \dots, a_\ell) = 0 \quad   \text{if one} \quad a_i \in \Cc(2n)
\end{gather*}
for all $C \in \Cc(2n)$.  Since $\Cc(2n)$ commutes with $A$, such a
cochain is completely determined by its restriction
$\widetilde{\Omega} : A^\ell \to \Cc(2n) \otimes A$ verifying
\[C  \widetilde{\Omega}(a_1, \dots, a_\ell) =
\widetilde{\Omega}(a_1, \dots, a_\ell)  C.\]

for all $C \in \Cc(2n)$. It results that $\widetilde{\Omega}$ is
actually $A$-valued. Then the map $\Omega \rightsquigarrow
\widetilde{\Omega}$ induces an isomorphism \cite{GS}:
\[H^\ell(\Cc(2n) \otimes A) \simeq H^\ell(A).\]

To obtain the desired cohomology, that is, $H^\ell(\Cc(2n,2k))$ or
$H^\ell(\Cc(2n+1,2k))$, we use the isomorphism $\phi : \Cc(2n) \otimes
A \to \Cc(2n,2k)$ or $\Cc(2n+1,2k)$ in the Periodicity Lemma 1
(Lemma~\ref{3.1}): given a cochain $\Omega$ of $\Cc(2n) \otimes A$, we
introduce a cochain $\phi^*(\Omega)$ of $\Cc(2n,2k)$ or $\Cc(2n+1,2k)$
def\/ined~by
\[\phi^*(\Omega) (x_1, \dots, x_\ell) = \phi\big(\Omega\big(\phi^{-1} (x_1),
\dots,\phi^{-1} (x_\ell)\big)\big), \] for all $x_1, \dots, x_\ell \in
\Cc(2n,2k)$ or $\Cc(2n+1,2k)$. Then the map $\Omega \rightsquigarrow
\phi^*(\Omega)$ induces a cohomology isomorphism.

It remains to compute the cohomology of $\Cc(1,2k) = \Cc(1)
\otimes_{\ZZ_2} \Wb_{2k}$. Since $\Cc(1,2k) = \Sk_2 \ltimes \Wb_{2k}$,
this is a particular case of a result in~\cite{Alev} where the
cohomology of $G \ltimes \Wb_{2k}$ is given for $G$ a~f\/inite group of
symplectic linear transformations. There is an improved version of
this result in~\cite{P07}, that allows a better management of
cocycles. Denote by $P$ the generator of $\Cc(1)$ satisfying $P^2 =
1$. One has:

\begin{prop}[\cite{Alev, P07}]\label{5.2}\qquad
\begin{enumerate}\itemsep=0pt

\item[$1.$] If $\ell >0$ and $\ell \neq 2k$, then
\[H^\ell (\Cc(1,2k)) =   \{0\}.
\]

\item[$2.$] $\dim\left(H^{2k}(\Cc(1,2k))\right) = 1.$ Moreover, there exists
  a normalized $\Cc(1)$-relative cocycle $\theta$ such that $H^{2k}
  (\Cc(1,2k)) = \CC \theta$ and
\[\theta (X_1, \dots, X_{2k}) = s^k(X_1 \wedge \dots \wedge X_{2k}) P, \qquad \text{for} \quad X_1, \dots, X_{2k} \in V_\un, \]
where $s$ is the canonical symplectic form on $V_\un$. If a
$2k$-cocycle $\theta'$ verifies $\theta'(X_1, \dots, X_{2k}) =
\theta(X_1, \dots, X_{2k}) $ for all $X_i \in V$, then $\theta' =
\theta \ {\rm{ mod }} \ B^{2k}$.
\end{enumerate}
\end{prop}

\begin{proof}
  See \cite{Alev} for the dimension of $H^\ell(\Cc(1,2k))$. See
  \cite{P07} for the last claims.
\end{proof}

\begin{cor} \label{5.3} \quad
\begin{enumerate}\itemsep=0pt

\item[$1.$] If $\ell >0$ and $\ell \neq 2k$, then
\[H^\ell (\Cc(2n+1,2k)) =  \{0\}.\]

\item[$2.$] Denote by $\omega_1, \dots, \omega_{2n}, P$ the canonical
  generators of $\Cc(2n+1)$ realized in $\Cc(2n+1,2k)$. Then there
  exists a cocycle $\Omega$ such that
such that $H^{2k} (\Cc(2n+1,2k)) =
  \CC \Omega$,
\begin{gather*}
    \Omega (X_1, \dots, X_{2k}) = i^n s^k(X_1 \wedge \dots \wedge
  X_{2k}) \omega_1 \dots \omega_{2k+1}, \qquad \text{for}\quad X_i \in
  V_\un, \\    \Omega (X_1, \dots, X_{2k}) = 0, \qquad \text{if one}\quad  X_i \in
  V_\ze.
\end{gather*}

If a $2k$-cocycle $\Omega'$ verifies $\Omega'(X_1, \dots, X_{2k}) =
\Omega(X_1, \dots, X_{2k}) $ for all $X_i \in V$, then $\Omega' =
\Omega$ ${\rm{mod}} \ B^{2k}$.
\end{enumerate}
\end{cor}

\begin{proof}
  Proposition \ref{5.2} provides a cocycle $\theta$ that allows us to
  construct a cocycle $\widehat{\theta}$ of $\Cc(2n) \otimes
  \Cc(1,2k)$ such that:
\[\widehat{\theta} (C_1 \otimes x_1, \dots, C_{2k} \otimes x_{2k}) = C_1  \cdots C_{2k} \otimes \theta(x_1, \dots, x_{2k})\]
for $x_1, \dots, x_{2k} \in \Wb_{2k}$, $C_1, \dots, C_{2k} \in
\Cc(2n)$. Next we compute $\Omega = \phi^*(\widehat{\theta})$ using
formulas in the proof of Lemma~\ref{3.1}:
\begin{gather*}
\Omega(X_1, \dots, X_{2k})  =  \phi\big( \widehat{\theta}(i^n \omega_1
 \cdots  \omega_{2n} \otimes X_1, \dots,i^n \omega_1
\cdots  \omega_{2n} \otimes X_{2k} )\big) \\
\phantom{\Omega(X_1, \dots, X_{2k})}{} =  \phi\big( (i^n)^{2k}
(\omega_1  \cdots  \omega_{2n})^{2k}  s^k(X_1 \wedge
\cdots \wedge X_{2k}) P\big)
\end{gather*}
for $X_1, \dots, X_{2k} \in V_\un$. Since $ (\omega_1 \cdots
\omega_{2n})^{2} = (-1)^n$ (see Section~\ref{Section00}), then
\begin{gather*}
  \Omega(X_1, \dots, X_{2k}) = \phi\big(s^k(X_1 \wedge \cdots \wedge X_{2k})
  P\big) = i^n s^k(X_1 \wedge \cdots \wedge X_{2k}) \omega_1  \cdots
  \omega_{2n}  P.\tag*{\qed}
\end{gather*}
\renewcommand{\qed}{}
\end{proof}

\begin{cor}
The Clifford--Weyl algebra $\Cc(2n+1,2k)$ is rigid if $k \neq 1$.
\end{cor}

Since $\dim \left( H^2(\Cc(2n+1,2)) \right) = 1$ and
$H^3(\Cc(2n+1,2))= \{0\}$, then $\Cc(2n+1,2)$ can be non trivially
deformed by a universal deformation formula (see Appendix~\ref{appendix2}). For
$\Cc(1,2)$, this formula is a particular case of a symplectic
ref\/lection algebras (see \cite{Etingof-Ginzburg}):

\begin{prop}
  Let $\Ac_\Lambda$ be the $\CC[\Lambda]$-algebra generated by $V_\ze
  = \CC \ P$ and $V_\un = \spa \{ E_+, E_- \}$ with relations:
  \[ [E_+,E_-]_\bl = -\tfrac14 + \Lambda P, \qquad P^2 =1 \qquad \text{and}\qquad
  PE_\pm = -E_\pm P.\]

  Then $\Ac_\Lambda$ is a non-trivial polynomial deformation of
  $\Cc(1,2)$ and a universal deformation formula.
\end{prop}

\begin{proof}
See \cite{Etingof-Ginzburg} or \cite {P07}.
\end{proof}

\section[Universal deformation formula of ${\mathcal C}(2n+1,2)$]{Universal deformation formula of $\boldsymbol{\Cc(2n+1,2)}$}\label{section6}

\begin{defn} \label{petal} Let $\apln$, $\lambda \in \CC$ be the
  algebra with generators $\omega_1, \dots, \omega_{2n+1}, E_\pm$ and
  relations:
\begin{gather*}
[E_+, E_-]_\bl = -\tfrac14 + i^n \lambda \omega_1 \cdots
  \omega_{2n+1},\\ 
   \omega_j \omega_k + \omega_k
  \omega_j = 2 \delta_{jk} \qquad (1 \leq j,k \leq 2n+1), \\
  E_\pm \omega_j = - \omega_j E_\pm \qquad (1 \leq j \leq 2n+1).
\end{gather*}
\end{defn}

\begin{defn} \label{petalg} The algebra $\agln$, when $\Lambda$ is a
  formal parameter, is def\/ined in a similar way: it is the algebra
  with generators $\omega_1, \dots, \omega_{2n+1}, E_\pm,\Lambda$ with
  $\Lambda$ central and same relations as $\apln$ with $\lambda$
  replaced by $\Lambda$. Note that $\agl(0) = \agl$.
\end{defn}

\subsection {\bf Construction using Ore extensions}

\begin{defn}
  Suppose that $R$ is an algebra, $\sigma$ an automorphism of $R$ and
  $\delta$ a $\sigma$-derivation of $R$, that is, a linear map $\delta
  : R \to R$ such that
\[ \delta(rs) = \delta(r) s +  \sigma(r) \delta(s) \]
for all $r,s \in R$.  Then the {\em Ore extension} $R[t]$ is the free
left $R$-module on the set $\{t^n \,|\, n \geq 0\}$, with
multiplication def\/ined by
\[ tr = \sigma (r)t + \delta(r) .\]
\end{defn}

Let $\Cc=\Cc(2n+1)$ be the Clif\/ford algebra in $2n+1$ generators,
$\omega_1, \dots, \omega_{2n+1}$. Consider the polynomial ring
$\Cc[\Lambda]$ where $\Lambda$ commutes with all elements of
$\Cc$. Elements of $\Cc[\Lambda]$ are denoted by $C(\Lambda)$.

Let $\tau$ be the automorphism of $\Cc[\Lambda]$ def\/ined by
\[\tau(\omega_r) = -\omega_r, \quad \forall \, r \qquad \text{and} \qquad \tau(\Lambda) = \Lambda.\]

The free $\Cc[\Lambda]$-module $\Cc[\Lambda][E_+]$ with basis $\{E_+^n
\,|\, n \in \NN\}$ gives us a a f\/irst Ore extension with
\[ E_+  C(\Lambda) = \tau(C(\Lambda)) E_+, \qquad \forall \, C(\Lambda) \in \Cc[\Lambda].\]

The following lemma is easy:

\begin{lem}\label{lem1}
There exists an automorphism $\sigma$ of the Ore extension
$\Cc[\Lambda][E_+]$ satisfying:
\[ \sigma(E_+) = E_+, \qquad \sigma(\omega_r) = - \omega_r, \quad \forall \, r \qquad \text{and} \qquad \sigma(\Lambda) = \Lambda. \]
\end{lem}

Let $\theta$ be the element $i^n \omega_1 \cdots \omega_{2n+1}
\Lambda$ in $\Cc[\Lambda]$. So $\theta$ commutes with $\Lambda$ and
$\omega_r$, $\forall \, r$ and anti-commutes with $E_+$.

Let $\Delta$ be the operator of $\Cc[E_+]$ def\/ined by
\[ \Delta(f) = \frac{f(E_+) - f(-E_+)}{2 E_+}, \qquad \forall \, f \in
\Cc[E_+]\] and $D$ be the operator of $\Cc[\Lambda][E_+]$ def\/ined by
\[
D(f(E_+)   C(\Lambda)) = \left(\frac14 \frac{df}{dE_+} - \Delta(f)
  \theta \right)   C(\Lambda), \qquad \forall \, f \in \Cc[E_+], \quad
C(\Lambda) \in \Cc[\Lambda].\]

\begin{lem}\label{lem2}
  One has $D(A B) = \sigma(A) D(B) + D(A) B$ for all $A, B \in
  \Cc[\Lambda][E_+]$.
\end{lem}

\begin{proof}
This is a straightforward verif\/ication.
\end{proof}

From Lemmas \ref{lem1} and \ref{lem2}, we can now construct a second
Ore extension  $\Cc[\Lambda][E_+][E_-]$ sa\-tisfying
\[ E_- A = \sigma (A) E_- + D(A), \qquad \forall \, A \in \Cc[\Lambda][E_+].
\]

It follows that:
\begin{gather}
  \notag  [E_+, E_-]_\bl = -\tfrac14 + \theta,\\ \label{rel1}
  \omega_k \omega_j + \omega_j \omega_k = 2 \delta_{jk} \qquad (1 \leq
  j,k \leq 2n+1),
   \\ \notag   E_\pm \omega_j = - \omega_j E_\pm \qquad (1
  \leq j \leq 2n+1).
\end{gather}

\begin{prop} \qquad
\begin{enumerate}\itemsep=0pt

\item[$1.$] The Ore extension $\Cc[\Lambda][E_+][E_-]$ and $\agln$ are
  isomorphic algebras.

\item[$2.$] A basis of $\agln$ is given by:
\[\left\{\omega^I   E_+^\alpha   E_-^\beta   \Lambda^r \, |\,  I \in \{0,1 \}^{2n+1}, \, \alpha, \beta, r \in \NN \right\},
\]
where $\omega^I = \omega_1^{i_1} \cdots \omega_{2n+1}^{i_{2n+1}}$ for
all $I = (i_1, \dots, i_{2n+1}) \in \{0,1 \}^{2n+1}$.
\end{enumerate}
\end{prop}

If $\Lambda$ is replaced by a small $\lambda$ ($\lambda \in \CC$) in
the def\/inition of $\agln$, the same procedure works to construct an
Ore extension of $\mathcal{C}[E_+][E_-]$, isomorphic to $\apln$. So

\begin{prop}
  A basis of $\apln$ is given by:
\[
\left\{\omega^I  E_+^\alpha  E_-^\beta \,|\, I \in \{0,1 \}^{2n+1}, \, \alpha, \beta \in \NN \right\}.
\]
\end{prop}

The algebra $\apln$ is the quotient $\agln/ I_\lambda$ where
$I_\lambda$ is the ideal $\agln \left( \Lambda - \lambda \right)$. As
a~particular case, setting $p = 2E_-$ and $q=2E_+$, we obtain:
\[\azn \simeq \Cc(2n+1,2) \simeq \agln/\Lambda \agln.\] Since $\agln =
\azn[\Lambda]$ as vector spaces, we obtain:

\begin{prop}\label{624}
  The algebra $\agln$ is a non-trivial polynomial deformation of the
  Clifford--Weyl algebra $\azn= \Cc(2n+1,2)$.
\end{prop}

\begin{proof}
  We just have to show that the deformation is non-trivial, but that
  results from the fact that the deformation cocycle is non-trivial by
  Corollary \ref{5.3}.
\end{proof}

\begin{rem}
  From Corollary \ref{5.3} and Lemma \ref{A.2.1}, this polynomial
  deformation $\agln$ is a~universal deformation formula of
  $\Cc(2n+1,2)$.
\end{rem}

\begin{cor}
  The center of $\agln$ is $\CC[\Lambda]$. Moreover, $\agln$ and
  $\apln$ are Noetherian algebras.
\end{cor}

\begin{proof}
  We have $\azn \simeq \mathcal{C}(2n+1,2)$ with center $\CC$ (Corollary
  \ref{2.5}). Let $\widetilde{a}$ be a central element of $\agln$. By
  Proposition \ref{624}, we can write $\widetilde{a} = a_0 + \Lambda
\widetilde{b}$ with $a_0 \in \azn$ and $\widetilde{b} \in
  \agln$. Therefore in $\agln$:
\[ x a_0 + \Lambda x \widetilde{b} = a_0 x + \Lambda \widetilde{b} x, \qquad \forall \, x \in \azn.\]

But $x a_0 = x \times a_0 + \Lambda \widetilde{c}$ and $a_0 x = a_0
\times x + \Lambda \widetilde{d}$ where $\times$ denotes the product
of $\azn$. So $a_0$ is central in $\azn$, henceforth $a_0 \in \CC$. It
follows $\tilde{b}$ is central in $\agln$ and repeating the same
argument, we obtain $\tilde{a} \in \CC[\Lambda]$. Finally, $\agln$ and
$\apln$ are Noetherian since they are constructed by Ore extensions of
Noetherian algebras~\cite{McCo}.
\end{proof}

In the sequel, we denote $\apl$ the algebra $\apl(0)$. The periodicity
of Clif\/ford algebras can be extended to the algebras $\agln$ and
$\apln$:

\begin{lem}[Periodicity Lemma~2] \label{isos}

  One has
\begin{gather*}
\agln \simeq \Cc(2n) \otimes \agl \simeq \Mcc_{2^n} \left(
  \agl \right) \qquad \text{and} \qquad
  \apln \simeq \Cc(2n) \otimes \apl \simeq \Mcc_{2^n} \left( \apl \right).
\end{gather*}
\end{lem}

\begin{proof}
  We denote by $P$ and $E_\pm$ the generators of $\agl(0)$ satisfying
  $P E_\pm = - E_\pm P$, $P^2 = 1$ and $[E_+, E_-]_\bl = - \frac14 +
  \Lambda P$. Let $\omega_1, \dots, \omega_{2n}$ be the generators of
  $\Cc(2n)$.

  We def\/ine $\omega_1', \dots, \omega_{2n+1}'$ and $E_\pm'$ elements
  of $\Cc(2n) \otimes \agl(0)$ by:
\begin{gather*}
 \omega_i' = \omega_i \otimes P \quad (1 \leq i \leq 2n),\qquad \omega_{2n+1}' = i^n
  \omega_1 \cdots \omega_{2n} \otimes P, \qquad E_\pm' = 1 \otimes E_\pm.
\end{gather*}

Using $(\omega_1 \cdots \omega_{2n})^2 = (-1)^n$, we check that
$\omega_1' , \dots, \omega_{2n+1}'$ verify the def\/ining relations of
$\Cc(2n+1)$ and anti-commute with $E_\pm'$. The relation $[E_+',
E_-']_\bl = -\frac14 + i^n \Lambda   \omega'_1 \cdots \omega'_{2n+1}$
results from $1 \otimes P = i^n \omega_1' \cdots \omega_{2n+1}'$.

Finally, this last equality and the fact that $\omega_i = \omega_i' (1
\otimes P)$ imply that $\omega_i'$ ($1 \leq i \leq 2n+1$) and $E_\pm'$
generate the algebra $\Cc(2n) \otimes \agln$.

On the other hand, if $\omega_1, \dots, \omega_{2n+1}$, $E_\pm$ are
the generators of $\agln$, we def\/ine $\omega_1', \dots, \omega_{2n}'$,
$E_\pm'$ and $P'$ by:
\begin{gather*}
 P' = i^n \omega_1 \cdots \omega_{2n+1}, \qquad \omega_{i}' = \omega_i P' \quad (1 \leq i \leq 2n),
\qquad E_\pm' = E_\pm.
\end{gather*}

Since $P'$ commutes with $\omega_i$, it commutes with
$\omega_i'$. Since $E_\pm'$ anti-commute with $\omega_i$, they
anti-commute with $P'$ and commute with $\omega_i'$. The equality
$P'^2 = 1$ follows from $(\omega_1 \cdots \omega_{2n+1})^2 = (-1)^n$
and we conclude $\omega_i'^2 = 1$. Moreover $\omega_i'$ anti-commutes
with $\omega_j'$ for $i \neq j$ and $[E_+', E_-']_\bl = -\frac14 + i^n
\Lambda   \omega_1 \cdots \omega_{2n+1} = -\frac14 + \Lambda P'$.

All def\/ining relations of $\Cc(2n) \otimes \agl(0)$ are
satisf\/ied. Moreover $\omega_i \!= \!\omega_i' P'$ and $i^n (\omega_1'
\cdots \omega_{2n+1}') P'$ $ = i^{2n} (\omega_1 \cdots \omega_{2n})^2
P'^{2n} \omega_{2n+1} = \omega_{2n+1}$. So we conclude that
$\omega_i'$, $E_\pm$ and $P'$ generate $\agln$.

This ends the proof that $\agln \simeq \Cc(2n) \otimes \agl$. Since
$\Cc(2n) \simeq \Mcc_{2^n} (\CC)$, then $\agln \simeq \Mcc_{2^n}
\left( \agl \right)$.

A similar reasoning works for $\apln$.
\end{proof}

\begin{rem}
  The f\/irst isomorphism in Lemma \ref{isos} is not a surprise: if $A$
  is an algebra, all deformations of $\Mcc_k(A)$ are of type
  $\Mcc_k(A_\Lambda)$ where $A_\Lambda$ is a deformation of $A$ (see
  Appendix~\ref{appendix3}). Here, $\Cc(2n+1,2) \simeq \Mcc_{2^n} \left( \Cc(1,2)
  \right)$ and $\agl$ is a deformation of $\Cc(1,2)$.
\end{rem}

\section[Algebras ${\mathcal A}_\Lambda(n)$ and their representations]{Algebras $\boldsymbol{\agln}$ and their representations}\label{section7}

Let $\{E_+,E_-, Y, F,G\}$ be the usual generators of the Lie
superalgebra $\osp(1,2)$: one has $\osp(1,2)_\ze$ $= \spa\{ Y,F,G\}$,
$\osp(1,2)_\un = \spa \{ E_+,E_- \}$ and the commutation relations
\begin{gather*}
   [Y, E_\pm] = \pm \tfrac12 E_\pm, \qquad [Y, F] = F, \qquad  [Y, G] = -G, \qquad [F,G] = 2 Y, \\
  [F, E_+ ] = [G, E_-] = 0, \qquad [F, E_-] = - E_+, \qquad [G, E_+] = - E_-,\\
  [E_+, E_+] = F, \qquad [E_-, E_-] = -G, \qquad [E_+, E_-] = Y,
\end{gather*}
where $[\cdot, \cdot]$ denotes the super bracket.

Let $\Uc:=\Uc(\osp(1,2))$ be the enveloping algebra of
$\osp(1,2)$. Denote by $\theta \in \Uc$ the {\em ghost}:
\[ \theta:=\tfrac14 + [E_+, E_-]_\bl,
\]
where $[\cdot, \cdot]_\bl$ denotes the Lie bracket.

\begin{lem}[\cite{P90,ABP,ABF}]\label{eth}
The relation $\theta E_\pm = - E_\pm \theta$ holds in $\Uc$.
\end{lem}

\begin{proof}
We have $\theta = \frac14 + E_+ E_-- E_- E_+$, hence
\begin{gather*}
E_+ \theta = \tfrac14 E_+ + E_+^2 E_- - E_+ E_- E_+,\qquad
\theta E_+ = \tfrac14 E_+ + E_+ E_- E_+ - E_- E_+^2.
\end{gather*}
Therefore $E_+ \theta + \theta E_+ = \frac12 E_+ - [Y, E_+] =
0$. Similarly, we can prove that $E_- \theta = - E_- \theta$.
\end{proof}

Let us now consider the $\CC$-algebra $\ut$ def\/ined by:
\[ \ut := \left\langle \Es_+, \Es_-, \vartheta \mid [\Es_+, \Es_-]_\bl =
-\tfrac14 + \vartheta, \Es_\pm \vartheta = - \vartheta \Es_\pm
\right\rangle. \]

By Lemma \ref{eth}, the enveloping algebra $\Uc$ is a quotient of
$\ut$.

\begin{prop}[\cite{LP}]
There exists an algebra isomorphism between $\Uc$ and $\ut$.
\end{prop}

\begin{proof}
  Consider the subspace $V= V_\ze \oplus V_\un$ of $\ut$, with $V_\ze
  = \{ 0\}$ and  $V_\un = \spa \{ \Es_+, \Es_- \}$. Def\/ine a
  supersymmetric bilinear form $(\cdot, \cdot)$ on $V$
  (hence symplectic on $V_\un$) by:
\[ (\Es_+, \Es_-) = - \tfrac14,  \qquad (\Es_+, \Es_+) = (\Es_-, \Es_-) = 0.\]

The algebra $\ut$ is $\ZZ_2$-graded by the $\ZZ_2$-gradation of $V$.
Starting from $[\Es_+, \Es_-] = 2 \Es_+ \Es_- + \frac14 - \theta$ with
$\Es_\pm \theta = - \theta \Es_\pm$, we have:
\[[[\Es_+, \Es_-],\Es_\pm] = \pm \tfrac12 \Es_\pm. \]
Using the Jacobi identity, we get $[[\Es_+, \Es_+],\Es_-] = - 2
[[\Es_+,\Es_-],\Es_+] = - \Es_+$ and $[[\Es_+, \Es_+],\Es_+] =0$.

In the same way, $[[\Es_-, \Es_-],\Es_+] =\Es_-$ and $[[\Es_-,
\Es_-],\Es_-] =0$. We conclude that
\[ [[X, Y], Z] = 2  \left( (Y,Z) X + (X, Z) Y \right), \qquad \forall \, X, Y \in V_\un.\]

By Corollary~\ref{4.2}, we deduce a surjective algebra homomorphism
from $\Uc$ to $\ut$ and using Lemma~\ref{eth}, we f\/inish the proof.
\end{proof}

\begin{prop} \label{7.4}  \qquad
\begin{enumerate}\itemsep=0pt
\item[$1.$] $\agln \simeq \Cc(2n+1) \zdtimes \Uc$.

\item[$2.$] $\agln \simeq \Cc(2n) \otimes (\Sk_2 \ltimes \Uc) \simeq
  \Mcc_{2^n} (\Sk_2 \ltimes \Uc)$.

\end{enumerate}

\end{prop}

\begin{proof}
1.   Let $\omega_1, \dots, \omega_{2n+1}$ be the generators of
  $\Cc(2n+1)$. Here $\Cc(2n+1)$ is $\ZZ_2$-graded by $\dzd(\omega_i) =
  1$, $\forall\, i$. Def\/ine $\Lambda \in \Cc(2n+1) \zdtimes \Uc$ by
\[\Lambda = i^n \omega_1 \cdots \omega_{2n+1} \theta. \]

We see immediately that $\Lambda$ is a central element and that
$\omega_1, \dots, \omega_{2n+1}$, $E_\pm$ and $\Lambda$ satisfy the
def\/ining relations of $\agln$. Moreover, they generate $\Cc(2n+1)
\zdtimes \Uc$ since $\theta = i^n \omega_1 \cdots \omega_{2n+1}
\Lambda$. Then there exists a surjective algebra homomorphism from
$\agln$ to $\Cc(2n+1) \zdtimes \Uc$.

To def\/ine the inverse map, we introduce an element $\theta \in \agln$
by
\[\theta = i^n \omega_1 \cdots \omega_{2n+1} \Lambda. \]

To f\/inish the proof, we notice that elements $E_+$, $E_-$ and $\theta$
verify the def\/ining relations of $\ut \simeq \Uc$, hence $\omega_1,
\dots, \omega_{2n+1}$, $E_\pm$ and $\theta$ satisfy the def\/ining
relations of $\Cc(2n+1) \zdtimes \Uc$ and they generate $\agln$.

2. The parity of $\Uc$ is used to def\/ine the smash product $\Sk_2
  \ltimes \Uc = \Cc(1) \zdtimes \Uc$. To prove~(2), apply Lemma~\ref{isos} and~(1).
\end{proof}

\begin{rem}
  The algebra $\agl$ is a deformation of $\Cc(1,2) = \Sk_2 \ltimes
  \Wb_2$. Besides $\agl = \Sk_2 \ltimes \Uc$. So here is a particular
  case where a deformation of a smash product remains a smash
  product. Moreover, representations of $\agl$ are merely graded
  representations of $\Uc$.
\end{rem}

\begin{defn}
An algebra $A$ is a {\em FCR algebra} if:
\begin{enumerate}\itemsep=0pt
\item Every f\/inite-dimensional representation of $A$ is completely
  reducible.

\item The intersection of all kernels of f\/inite-dimensional
  representations is $\{0\}$.
\end{enumerate}

\end{defn}

If $A$ is an algebra and $V$ is an $A$-module, there is a
corresponding structure of $\Mcc_n(A)$-module on $\Vc =V^n$, since
$n\times n$-matrices with coef\/f\/icients in $A$ act on $n$-column
vectors with coef\/f\/icients in $V$. The following well-known Lemma shows
that all $\Mcc_n(A)$-modules are of this type:

\begin{lem} \label{7.6} Given an $\Mcc_n(A)$-module $\Vc$, there
  exists a corresponding $A$-module $V$ such that $\Vc \simeq V^n$,
  and $\ker(\Vc) = \Mcc_n(\ker(V))$.
\end{lem}

\begin{proof}
  Let $E_{ij}$ be the elementary matrices in $\Mcc_n$. They satisfy
  $\Id = \sum\limits_{i=1}^n E_{ii}$ and $E_{ii} E_{jj} = \delta_{ij}
  E_{ii}$. So $\Vc = \oplus_{i=1}^n W_i$ where $W_i = E_{ii}\cdot
  \Vc$. Since $A$ and $\Mcc_n$ commute, $W_i$ are $A$-modules. Since
  $E_{ij} \cdot W_j = E_{ii} \cdot (E_{ij} \cdot W_j)$, $E_{ij}$ maps $W_j$ into
  $W_i$ and $E_{ij} W_k = \{ 0\}$ if $k \neq j$. Using $E_{ji} E_{ij}
  = E_{jj}$, it results that $W_i$ and $W_j$ are isomorphic
  $A$-modules, for all $i$, $j$. Let $V$ be any of the $A$-modules~$W_i$,
  it is easy to verify that $\Vc \simeq V^n$ as
  $\Mcc_n(A)$-modules. For the second claim, use the following easy
  statement: if $J$ is a two-sided ideal of $\Mcc_n(A)$, there exists
  a two-sided ideal $I$ of $A$ such that $J = \Mcc_n(I)$.
\end{proof}

\begin{rem}
  The above correspondence between $A$-modules and $\Mcc_n(A)$-modules
  respects isomorphisms, irreducibility and direct sum decompositions.
\end{rem}

\begin{cor} \label{7.7}
If $A$ is an FCR algebra, then $\Mcc_n(A)$ is also a FCR algebra.
\end{cor}

\begin{proof}
  Apply Lemma~\ref{7.6} and the remark above.
\end{proof}

\begin{cor}\label{fcr}  \qquad
\begin{enumerate}\itemsep=0pt
\item[$1.$] The algebra $\agln$ is a FCR algebra.

\item[$2.$] All finite-dimensional representations of $\apln$ are completely
  reducible.
\end{enumerate}
\end{cor}

\begin{proof}
1.\ By Periodicity Lemma~2, $\agln \!\simeq \!\Mcc_{2^n}\!(\agl)$ and by
  Proposition~\ref{7.4}, $\agl \!\!\simeq \!\Sk_2 {\ltimes}
  \Uc$. Representations of $\Sk_2 \ltimes \Uc$ are simply
  $\ZZ_2$-graded representations of $\Uc$ and by~\cite{DH},
  f\/inite-dimen\-sio\-nal ones are completely reducible. Moreover by
  \cite{Behr}, the intersection of all kernels of f\/inite-dimensional
  representations of $\Sk_2 \ltimes \Uc$ is $\{0\}$. So $\Sk_2 \ltimes
  \Uc$ is FCR, then apply Corollary~\ref{7.7}.

2. Use $\apln \simeq \agln / (\Lambda - \lambda) \agln$.
\end{proof}

\begin{rem}
  By Lemma \ref{7.6}, the representation theory of $\agln \simeq
  \Mcc_{2^n}(\agl)$ is reduced to the representation theory of $\agl
  \simeq \Sk_2 \ltimes \Uc$ and therefore to the $\ZZ_2$-graded
  representation theory of $\Uc$.
\end{rem}

\begin{ex}
  To study f\/inite-dimensional representations of $\agln$, it is enough
  to study irreducible ones by Corollary \ref{fcr}. Irreducible
  representations of $\osp(1,2)$ are well-known, they are all
  $\ZZ_2$-graded: given $h \in \frac12 \NN$, there exists an
  irreducible representation on a $(4h+1)$-dimensional space $V_h$ to
  which corresponds an irreducible representation of $\agln$ given by
  $2^n \times 2^n$-matrices with coef\/f\/icients in $\Sk_2 \ltimes \Uc$
  acting on $\Vc_h = V_h^{2^n}$, therefore of dimension
  $2^n(4h+1)$. Alternatively, since $\agln \simeq \Cc(2n) \otimes
  \agl$, this representation is the natural action of $\Cc(2n) \otimes
  \agl$ on $\Vc_h = \Phi_n \otimes V_h$ where $\Phi_n$ is the spin
  representation of $\Cc(2n)$ (see Section~\ref{Section00}). This describes all
  f\/inite-dimensional representations of $\agln$.
\end{ex}

\begin{ex}\label{7.12}
  We will now construct examples of simple $\apl$-modules from the
  metaplectic representation of $\Wb_2$. Let $V = \CC[z]$ as in
  \cite{LP}. We def\/ine the operator $\Delta$ of $V$:
\[ \Delta(h) = \frac1z(h(z) - h(-z)), \qquad \forall \, h \in V.\]

Denote by $P$ the parity operator of $V$. Def\/ine operators $\rho_\lambda^\pm$
by:
\[\rho_\lambda^+= \frac12 \frac{d}{dz} - \lambda \Delta,   \qquad
\rho_\lambda^- = -\frac12 z.\]

Then $[\rho_\lambda^+, \rho_\lambda^-]_\bl = - \frac14 + \lambda P$,
$\rho_\lambda^\pm P = -P \rho_\lambda^\pm$ and $P^2 = 1$. In this way,
we obtain a representation~$ \rho_\lambda$ of~$\apl$ in $V$ such that:
\[\rho_\lambda(E_\pm) = \rho_\lambda^\pm, \qquad \rho_\lambda(P) =P.\]

We recover exactly the $\ZZ_2$-graded $\osp(1,2)$-Verma module
$\Vc_{\lambda - \frac14}$ of highest weight $\lambda - \frac14$.

If $\lambda \neq h + \frac14$, $2h \in \NN$, then $\rho_\lambda(E_+)$
does not vanish and the corresponding module is simple. If $\lambda =
h + \frac14$, $2h \in \NN$, we have $\rho_\lambda(E_+)(z^{4h+1})
=0$. Therefore $W_h = \spa \{ z^\ell , \ell \geq 4h+1 \}$ is a simple
submodule of dominant weight $-\left(h+\frac12\right)$, the quotient
$V/W_h$ is the simple $\osp(1,2)$-module of dimension $4h+1$ and the
module $(V, \rho_{h+\frac14})$ is a non-trivial extension of $W_h$ by
$V/W_h$ (see \cite{LP} for more details).
\end{ex}

Denote by $V_\lambda$ the $\apl$-module just built. Using $\apln =
\Mcc_{2^n}(\apl)$, def\/ine a corresponding $\apln$-module by setting
$V_\lambda(n) = \Phi_{2n} \otimes V_\lambda$ where $\Phi_{2n}$ is the
spin representation of $\Cc(2n)$. When $\lambda \neq \left( h +
  \frac14 \right)$, $2h \in \NN$, we obtain a simple $\apln$-module.
When $\lambda = h + \frac14$, $2h \in \NN$, we obtain an
indecomposable $\apln$-module with a unique simple submodule and a
unique simple quotient of dimension $2^{n} (4h+1)$.

Since $\agln \simeq \agln / (\Lambda - \lambda) \agln$, these modules
are $\agln$-modules.

\begin{rem}
  When $A$ is a $\ZZ_2$-graded algebra, $\Mcc_n(A) = \Mcc_n \otimes A$
  has a natural $\ZZ_2$-gradation induced by the gradation of $A$ and
  $\deg(M) = \zero$ for all $M \in \Mcc_n$. But algebras $\Sk_2
  \ltimes \Mcc_n(A)$ and $\Mcc_n(\Sk_2 \ltimes A)$ have the same
  underlying vector space. It is easy to verify that they coincide as
  algebras. Using Proposition \ref{7.4}, Remark \ref{3.4} and Lemma
  \ref{7.6}, it results that all representations of $\agln$ are graded
  and obtained from graded representations of $\Uc$.
\end{rem}

\section[Algebras ${\mathcal A}_\alpha(n)$]{Algebras $\boldsymbol{\apln}$}\label{section8}

We keep the notation of last Section. Write $\g=\osp(1,2)$ as $\g=\gO
\oplus \gI$ where $\gO = \spa\{ Y,F,G \}$ and $\gI = \spa \{ E_+, E_-
\}$, $\Uc = \Uc(\g)$ its enveloping algebra and $\theta = \frac14 + [
E_+, E_-]_\bl$ the ghost. We have $\Zc(\g) = \CC[C]$ where $\Zc(\g)$
denotes the center of $\Uc$, $C = \theta^2 - \frac{1}{16}$ and
$\Zc(\gO) = \CC[Q]$ where $\Zc(\gO)$ denotes the center of $\Uc(\gO)$,
$Q = \left( \theta - \frac14 \right) \left( \theta + \frac34 \right)$
\cite{P90, ABP}. For $c \in \CC$, let $\Bc_c:= \Uc / (C-c)\Uc$.

Let us consider the $\CC$-algebra $\apl:=\apl(0)$. Recall that:
\[ \apl = \left\langle  E_+, E_-, P \,|\, P^2 = 1, [E_+, E_-]_\bl = -\tfrac14 + \lambda P, E_\pm P = - P E_\pm \right\rangle. \]

If $\lambda = 0$, $\Ac_0$ is the Clif\/ford--Weyl algebra $\Cc(1,2) =
\Sk_2 \ltimes \Wb_2$. In general:

\begin{prop}\label{ab}
  One has $\apl \simeq \Bc_{\lambda^2 - \frac1{16}}$ whenever $\lambda
  \neq 0$.
\end{prop}

\begin{proof}
  For $u \in \Uc$, we denote by $\overline{u}$ its class in
  $\Bc_{\lambda^2 - \frac1{16}}$. Therefore
  $[\overline{E_+},\overline{ E_-}]_\bl = -\frac14 +
  \overline{\theta}$ and $ \overline{E_\pm} \ \overline{\theta} = \pm
  \overline{\theta} \ \overline{E_\pm}$. Moreover,
  $\overline{C}=\lambda^2 - \frac1{16}= \overline{\theta}^2 -
  \frac{1}{16}$.  Setting $P = \frac1\lambda \overline{\theta}$, one
  recovers exactly the def\/ining relations of $\apl$ and a map from
  $\apl$ onto $\Bc_{\lambda^2 - \frac1{16}}$.

  For the inverse map, one can check that elements $E_+$ and $E_-$ in
  $\apl$ generate a superalgebra isomorphic to $\g$, hence a
  homomorphism $\rho$ from $\Uc$ to $\apl$. We have $\rho(\theta) =
  \lambda P$, so $\rho$ is surjective. Since $\rho \left(C - \lambda^2
    + \frac1{16} \right) = 0$, one can def\/ine the inverse map
  $\overline{\rho}$ from $\Bc_{\lambda^2 - \frac1{16}}$ onto $\apl$.
\end{proof}

The structure of the algebra $\apl$ is deduced from the Proposition
above and \cite{P90}.

\begin{prop} \label{8.2} \qquad

\begin{enumerate}\itemsep=0pt
\item[$1.$] If $\lambda^2 \neq \left(h + \frac14 \right)^2$, $2h \in \NN$,
  then $\apl$ is a simple algebra.

\item[$2.$] If $\lambda^2 = \left(h + \frac14 \right)^2$, $2h \in \NN$, then
  $\apl$ is a primitive algebra. Moreover, there exists a unique
  non-zero two-sided ideal $I_\lambda$ in $\apl$ of codimension
  $\left( 4h +1 \right)^2$, with $I_\lambda = \ker(V_h)$ and $V_h$ is
  the simple $\osp(1,2)$-module of dimension $4h +1$.
\end{enumerate}
\end{prop}

\begin{proof}
  It is proved in \cite{P90} that $\Bc_c$, $c \neq 0$ has the
  following properties:

\begin{itemize}\itemsep=0pt

\item if $c \neq h\frac{(2h+1)}{2}$, $2h \in \NN$, then $\Bc_c$ is
  $\ZZ_2$-simple;

\item if $c = h\frac{(2h+1)}{2}$, $2h \in \NN$, then $\Bc_c$ is
  primitive. Moreover, there exists a unique non-zero $\ZZ_2$-graded
  two-sided ideal with codimension $\left(4 h+1 \right)^2$ that is the
  kernel of the simple $\osp(1,2)$-module of dimension $4h+1$.

\end{itemize}

If $\lambda = 0$, then $\Ac_0 \simeq \Cc(1,2)$ is simple.

If $\lambda \neq 0$, then $\apl \simeq \Bc_{\lambda^2 -
  \frac1{16}}$. It is enough to show that any two-sided ideal of
$\Bc_{\lambda^2 - \frac1{16}}$ is $\ZZ_2$-graded and then translate
the results just above in term of $\lambda$. So, let $I$ be a
two-sided ideal of $\Bc_{\lambda^2 - \frac1{16}}$. We set $P =
\frac1\lambda \overline{\theta}$. We have $P^2 =1$ and $P b P =
(-1)^{\dzd(b)} b$, $\forall\, b \in \Bc_{\lambda^2 - \frac1{16}}$. If $a
= a_\ze + a_\un \in I$, it follows $PaP = a_\ze - a_\un \in I$,
therefore $a_\ze$ and $a_\un \in I$.
\end{proof}

\begin{cor} \label{8.3} \qquad

\begin{enumerate}\itemsep=0pt
\item[$1.$] If $\lambda^2 \neq \left(h + \frac14 \right)^2$, $2h \in \NN$,
  then $\apln$ is a simple algebra.

\item[$2.$] If $\lambda^2 = \left(h + \frac14 \right)^2$, $2h \in \NN$, then
  $\apln$ is a primitive algebra. Moreover, there exists a~unique
  non-zero two-sided ideal in $\apln$ of codimension $2^{2n} \left( 4h
    +1 \right)^2$, that is the kernel of the irreducible
  representation of dimension $2^{n} \left( 4h +1 \right)$.
\end{enumerate}
\end{cor}

\begin{proof}
  By Lemma \ref{isos}, $\apln \simeq \Mcc_{2^n} \left( \apl \right)$,
  so two-sided ideals of $\apln$ are all of type $\Mcc_{2^n}(I)$, $I$
  is a two-sided ideal of $\apl$. Then apply Proposition~\ref{8.2} and
  Example \ref{7.12}.
\end{proof}

\begin{rem}
  We have $\agl(n) \simeq \agl(n) / (\Lambda - \lambda)
  \agl(n)$. Moreover, the center of $\agl(n)$ is $\CC[\Lambda]$, so
  Corollary \ref{8.3} lists all primitive quotients of $\agl(n)$.
\end{rem}

The algebra $\Ac_0 = \Cc(1,2)$ is a quotient of $\Uc(\osp(2,2))$ (see
Proposition \ref{4.1}). More generally:

\begin{prop} \label{8.7}

  The algebra $\agl$ is a quotient of $\Sk_2 \ltimes
  \Uc(\osp(2,2))$. Moreover, the Casimir operator of $\osp(2,2)$ $($see
  {\rm \cite{ABP}}$)$ vanishes in this quotient.
\end{prop}

\begin{proof}
  Using the notation in Section~\ref{section6}, we consider in $\agl$:
\[ K = -\tfrac14 \omega_1 + \Lambda, \qquad H_\ze = \CC K, \qquad H_\un = \spa
\{ E_+, E_- \}, \qquad H = H_\ze \oplus H_\un.\]

Def\/ine $(\cdot| \cdot)$ a supersymmetric bilinear form on $H$ by
$(K|K) = \frac18$ and $(E_+ | E_-) = - \frac14$. It is easy to check
that relations (\ref{six}) hold in $H$, so by Corollary~\ref{4.2}, the
subalgebra of $\agl$ generated by $H$ is a quotient of
$\Uc(\osp(2,2))$. Now, the subalgebra of $\agl$ generated by $H$ and
$\omega_1$ is $\agl$ itself, and it is clearly a quotient of $\Sk_2
\ltimes \Uc(\osp(2,2))$. The second claim results from a direct
computation using the Casimir formula given in~\cite{ABP}.
\end{proof}

\begin{cor}
  Any graded $\osp(1,2)$-module can be extended to an
  $\osp(2,2)$-module $($with same underlying space$)$.
\end{cor}

\begin{proof}
  First, remark that given a graded algebra $A$, graded $A$-modules
  and $\Sk_2 \ltimes A$-modules are exactly the same notion.

  Now, start with a graded $\osp(1,2)$-module with parity $P$. Recall
  that $[E_+, E_-]_\bl = -\frac14 + \theta$ with $\theta E_\pm = -
  E_\pm \theta$. We def\/ine $\Lambda = \theta P$ and $\omega_1 = P$ to
  obtain a graded $\agl$-module. By Proposition~\ref{8.7}, this module
  is a $\Sk_2 \ltimes \Uc(\osp(2,2))$-module, therefore a graded
  $\osp(2,2)$-mo\-dule.
\end{proof}

\begin{rem}
  Let $\Cs$ be the Casimir element of $\Uc(\osp(2,2))$. It is proved
  in \cite{ABP} that a simple $\osp(2,2)$-module is still simple as an
  $\osp(1,2)$-module if, and only if, $\Cs = 0$.
  \end{rem}

\appendix
\section{Appendix}\label{appendix1}

For the convenience of the reader, we recall here some notions of
Hochschild cohomology theory relating it to Gerstenhaber deformation
theory of (associative) algebras \cite{G64,Gre}. See \cite{BFFLS, DS}
for applications of deformation theory to quantization.

Let $A$ be an (associative) algebra. By Hochschild cohomology of $A$,
we mean Hochschild cohomology with coef\/f\/icients in $A$, def\/ined as
follows.

For $k > 0$, $k$-{\em cochains} are $k$-linear maps from $A^k$ to
$A$. When $k=0$, $0$-cochains are simply elements of $A$. We denote by
$M^k(A)$ the space of $k$-cochains and by $M(A) = \oplus_{k \geq 0}
M^k(A)$, the space of cochains. We def\/ine the {\em Hochschild
  coboundary operator} $d$ acting on $M(A)$ by:
\begin{itemize}\itemsep=0pt
\item if $a \in A = M^0(A)$, $da = - \ad(a)$ where $\ad(a)(b):=[a,b]$,
  for all $a$,$b \in A$;
\item if $\Omega \in M^k(A)$, $k>0$:
 \begin{gather*}
   d \Omega(a_1, \dots,a_{k+1}) = a_1 \Omega(a_2, \dots, a_{k+1}) -
   \Omega(a_1 a_2, a_3, \dots, a_{k+1}) \\
   \phantom{d \Omega(a_1, \dots,a_{k+1}) =}{} +   \Omega(a_1, a_2 a_3, \dots, a_{k+1}) +\cdots +(-1)^{k+1} \Omega(a_1, \dots, a_k) a_{k+1}.
\end{gather*}
\end{itemize}

One has $d^2 =0$. Let $B^0(A)= \{0\}$, $B^k(A)= dM^{k-1} (A)$,
$k>0$. Set $Z^k(A) = \ker\big(d|_{M^k(A)} \big)$, $k \geq 0$ and
$H^k(A) = Z^k(A) / B^k(A)$. Elements of $B^k(A)$ (resp. $Z^k(A)$)
are $k$-{\em coboundaries} (resp.\ $k$-{\em cocycles}) and $H^k(A)$
is the $k^{\text{th}}$-space of Hochschild cohomology of $A$. Note
that $H^0(A)$ is the center of $A$. Here are some examples of
algebras $A$ such that $H^k(A) = \{0\}$ for all $k > 0$: semi-simple
algebras (e.g. algebras of f\/inite groups, algebras of complex
matrices, Clif\/ford algebras), Weyl algebras, etc.

A {\em deformation} of $A$ with formal parameter $\Lambda$ is a
$\CC[[\Lambda]]$-algebra structure on $A[[\Lambda]]$ def\/ined~by:
\begin{gather*}
a \sta b = ab + \sum_{n \geq 1} \Lambda^n \Omega_n(a,b), \quad \forall\,
a,b \in A, \qquad \Omega_n \in M^2(A), \quad \forall \, n.
\end{gather*}

The associativity of $\sta$ can be reinterpreted in terms of
Hochschild cohomology: $\Omega_1 \in Z^2(A)$ and when $\Omega_1 \in
B^2(A)$, it can be removed by an equivalence, i.e.\ an isomorphism of
$\CC[[\Lambda]]$-algebras. When $H^2(A) = \{0\}$, repeating the same
argument, it results that any deformation is equivalent to the initial
product, so $A$ is rigid. For instance, all algebras we just mentioned
above are rigid. Second, the conditions on $\Omega_n$, $n \geq 2$ can
be written in terms of 3-cohomology, and it results that if $H^3(A) =
\{0\}$, then given any $\Omega_1 \in Z^2(A)$, there exits a
deformation with leading cocycle $\Omega_1$.

These two results are known as the rigidity and integrability
theorems.

\section{Appendix}\label{appendix2}

The terminology and results presented in this Appendix are rather
standard, but for the sake of completeness we include them here with
proofs.

Let $A$ be an associative algebra with product $m_0$. Let $\Mb (A) =
\sum\limits_{k \geq 0} \Mb^k(A)$ be the space of multilinear maps from $A$ to
$A$. The space $\Mb(A)$ is graded, $\Mb^{(k)}:= \Mb^{k+1}(A)$ and
endowed with the Gerstenhaber bracket, it is a graded Lie algebra~\cite{PU07}. Let $d = - \ad(m_0)$. Since $d^2 = 0$, $d$ def\/ines a
complex on $\Mb(A)$, the {\em Hochschild cohomology complex} of $A$
(see~\cite{GS}). Let $\Zb^2(A)$ be the set of 2-cocycles, $B^2(A)$ the
2-coboundaries, and $H^2(A)$ chosen such that $\Zb^2(A) = B^2(A)
\oplus H^2(A)$.

Given two vector spaces $V$ and $W$, a {\em formal map} $F : V \to W$
is a power series $F = \sum\limits_{k \geq 0} F_k$ where $F_k$ is a
homogeneous polynomial function of degree $k$ from $V$ to $W$. In the
sequel, we will need essentially formal maps $F : H^2(A) \to \Mb(A)$
and we def\/ine a graded Lie algebra bracket coming from the one def\/ined
on $\Mb(A)$ by:
\[
[F, F'] = \sum_{k \geq 0} \sum_{r+s=k} [F_r, F_s'] \qquad \text{for} \qquad F = \sum_{k \geq 0}F_k,
 \qquad F' = \sum_{k \geq 0}F_k'
\]
with $[F_r, F_s'](h) = [F_r(h), F_s'(h)]$, $\forall\, h \in H^2(A)$.

\begin{defn}
  A {\em universal deformation formula} of $A$ is a formal map $F :
  \Zb^2(A) \to \Mb^2(A)$ such that:
\begin{enumerate}\itemsep=0pt
\item[$1)$] $F = m_0 + \Id_{H^2(A)} + \sum\limits_{k \geq 2} F_k$,

\item[$2)$] $[F,F]=0$.
\end{enumerate}
\end{defn}

If $F$ is a universal formula of deformation, $\lambda$ a formal
parameter and $h \in H^2(A)$, then $m_h^\lambda := F(\lambda h) = m_0
+ \lambda h + \sum\limits_{k \geq 2} \lambda^k F_k$ is a deformation of
$m_0$. More generally, if we have a formal curve in
$H^2(A)[[\lambda]]$, $\tilde{h} = \sum\limits_{n \geq 1} \lambda^n h_n$, then
\[ m_{\tilde{h}}^\lambda := F(\tilde{h}(\lambda)) = m_0 + \lambda h_1
+ \sum_{k \geq 2}  \lambda^k \underset{i_1, \dots, i_n \geq 1, 1 \leq
  n \leq k}{\sum_{i_1 + \cdots + i_n =k}} F_n (h_{i_1}, \dots,
h_{i_n}) \] is a deformation of $m_0$. The Lemma below is simply a
translation of the classical criterion of integrability:

\begin{lem}\label{A.2.1}
Let $\Db^2(A)$ be a complementary subspace of $\Zb^2(A)$ in
$\Mb^2(A)$. If $H^3(A) = \{ 0 \}$, then there exists a universal
deformation formula
\[F = m_0 + \Id_{H^2(A)} + \sum_{k \geq 2} F_k, \qquad \text{with} \quad F_k \in \Db^2(A), \quad \forall \, k \geq 2.\]
\end{lem}

\begin{proof}
  Let $\sigma$ be a section of $d : \Mb^2(A) \to \Bb^3(A)$ such that
  $\sigma \circ d$ is the projection onto $\Db^2(A)$ along
  $\Zb^2(A)$. Step by step, we construct $F$ verifying $[F,F] = 0$, $F
  = m_0 + \Id_{H^2(A)} + \sum\limits_{k \geq 2} F_k$: f\/irst, we f\/ind $d(F_2)
  = \frac12 [ \Id_{H^2(A)}, \Id_{H^2(A)}]$. Since $[ \Id_{H^2(A)},
  \Id_{H^2(A)}]$ is valued in $\Zb^3(A) = \Bb^3(A)$, def\/ine a suitable
  $F_ 2 = $  $ \frac12 \sigma \circ [ \Id_{H^2(A)},
  \Id_{H^2(A)}]$ (remark that $d \circ \sigma = \Id_{\Bb^3(A)}$). It
  is easy to see that the remaining $F_k$ can be constructed by the
  same procedure.
\end{proof}

{\samepage
\begin{lem} \label{A.2.2} \qquad
\begin{enumerate}\itemsep=0pt
\item[$1.$] Let $m^\lambda$ be a deformation of $m_0$. Then, up to
  equivalence, $m^\lambda$ can be written as:
\[m^\lambda = m_0 + h(\lambda) + d(\lambda), \qquad \text{with} \quad h(\lambda) \in \lambda H^2[[\lambda]], \quad d \in \lambda^2 \Db^2[[\lambda]].\]

\item[$2.$] If $m'^\lambda$ is another deformation with
\[m'^\lambda = m_0 + h(\lambda) + d'(\lambda), \qquad \text{with}  \quad d'(\lambda) \in \lambda^2 \Db^2[[\lambda]],\]
then $d'(\lambda) = d(\lambda)$.
\end{enumerate}
\end{lem}}

\begin{proof} 1. Up to equivalence, we can assume that the leading cocycle of
  $m^\lambda$ is in $H^2(A)$, $m^\lambda = m_0 + \lambda h_1 +
  \lambda^2 C_2+ \cdots$.

We have $C_2 = d_2 + h_2 + b_2$, $d_2 \in \Db^2(A)$, $h_2 \in H^2(A)$
and $b_2 \in \Bb^2(A)$. We can assume that $b_2 = 0$, therefore
$m^\lambda = m_0 + (\lambda h_1 + \lambda^2 h_2) + \lambda^2 d_2 +
\lambda^3 C_3+ \cdots$. Repeat the same argument to obtain the result.

2. Let $m^\lambda = m_0 + \lambda h_1 + \lambda^2 (h_2 + d_2) +
  \cdots$, $m_\lambda '= m_0 + \lambda h_1 + \lambda^2 (h_2 + d_2') +
  \cdots$, then $d(h_2 + d_2) = \frac12 [h_1, h_1] = d(h_2 + d_2')$,
  hence $d(d_2) = d(d_2')$ and that implies $d_2 - d_2' \in \Zb^2(A)
  \cap \Db^2(A) = \{0\}$. Apply repeatedly the same reasoning to
  obtain $m_\lambda' = m^\lambda$.
\end{proof}

\begin{prop}
  Assume that $H^3(A) = \{0\}$. Let $F$ be a universal deformation
  formula and~$m^\lambda$ a deformation. Up to equivalence, there
  exists a formal curve $h(\lambda)$ in $H^2(A)[[\lambda]]$ such that
  $h(0) = 0$ and $m^\lambda = F(h(\lambda))$. In other words, $F$
  characterizes all deformations of $m_0$ up to equivalence and up to
  change of formal parameter.
\end{prop}

\begin{proof}
  The existence is given by the Lemma \ref{A.2.1}. Up to equivalence,
  we can assume that $m^\lambda= m_0 + h(\lambda) + d(\lambda)$ where
  $h(\lambda) \in \lambda H^2[[\lambda]]$ and $d(\lambda) \in
  \lambda^2 H^2[[\lambda]]$ (Lemma \ref{A.2.2}). But the deformation
  $m'^\lambda = F(h(\lambda))$ can be written as $m'^\lambda = m_0 +
  h(\lambda ) + d'(\lambda)$ with $d'(\lambda) \in
  \lambda^2H^2[[\lambda]]$. Henceforth $m'^\lambda = m^\lambda$ by
  Lemma~\ref{A.2.2}.
\end{proof}

\section{Appendix}\label{appendix3}

Let $A$ be an algebra and $\Agl$ be a deformation of $A$ with product
$\sta$. The underlying space of $\Agl$ is $A [[\Lambda]]$, and it is
easy to check that $\Mcc_n(A[[\Lambda]]) = \Mcc_n(A)[[\Lambda]]$. Then
$\Mcc_n(\Agl)$ is a deformation of $\Mcc_n(A)$, the product is the
natural one, def\/ined by $(aM)\sta (a'M')= (a \sta a') MM'$, $\forall\,
a, a' \in A$, $M, M' \in \Mcc_n(A)$. Conversely:

\begin{prop}
  Any deformation of $\Mcc_n(A)$ is equivalent to a deformation
  $\Mcc_n(\Agl)$ with~$\Agl$ a deformation of $A$.
\end{prop}

This result is known, but since we have not been able to f\/ind a
reference, we give a short proof.

\begin{proof}
  We refer to \cite{GS} for relative deformation theory with respect
  to a separable subalgebra. In the present case, the separable
  subalgebra of $\Mcc_n(A)$ is $\Mcc_n$, and any deformation is
  equivalent to a deformation with normalized $\Mcc_n$-relative
  cochains \cite{GS}, that is, cochains $\Omega : \left( \Mcc_n(A)
  \right)^2 \to \Mcc_n(A)$ that verify for all $M \in \Mcc_n, a_1, a_2
  \in A$:
\begin{gather*}
 \Omega(M a_1, a_2) = M \Omega(a_1, a_2), \qquad \Omega(a_1 M, a_2) =
  \Omega(a_1, M a_2), \\  \Omega(a_1, a_2M) = \Omega(a_1, a_2)M, \qquad \text{and}
   \qquad \Omega(x_1, x_2) = 0 \qquad \text{if one} \quad x_i \in \Mcc_n.
\end{gather*}
Since $\Mcc_n$ and $A$ commute, such a cochain is completely
determined by its restriction $\widetilde{\Omega} : A^2 \to \Mcc_n(A)$
that verif\/ies $M \widetilde{\Omega} (a_1, a_2) = \widetilde{\Omega}
(a_1, a_2)M$, $\forall \, M \in \Mcc_n$, $a_1, a_2 \in A$, and is
therefore $A$-valued. Summarizing, up to equivalence, we have a new
product $\sta$ that satisf\/ies
\begin{gather*}
  M_1 \sta M_2 = M_1 M_2,  \qquad\! M_1 \sta a = a \sta M_1 = a M_1, \qquad\! (a_1 M_1) \sta (a_2 M_2) = (a_1 \! \sta a_2)M_1M_2
\end{gather*}
for all $M_1, M_2 \in \Mcc_n$, $a, a_1, a_2 \in A$ and
\begin{gather*}
 a_1 \sta a_2 = a_1 a_2 + \sum_{n\geq 1} \Lambda^n C_n(a_1, a_2), \qquad
\forall \, a_1, a_2 \in A
\end{gather*} with $C_n : A^2 \to A$. So $\sta$ def\/ines a
deformation $\Agl$ of $A$.

Now, we will prove that our initial deformation $\sta$ of $\Mcc_n(A)$
is exactly the deforma\-tion~$\Mcc_n(\Agl)$: it is enough to show that
$(aM) \sta (a'M')$ is the product of~$(aM)$ and $(a'M')$ in~$\Mcc_n(\Agl)$, for all $a, a' \in A$, $M,M' \in \Mcc_n$. But this is
true since $aM \sta a'M' = (a\sta a')MM'$, that is exactly the product
of $\Mcc_n(\Agl)$.
\end{proof}

\pdfbookmark[1]{References}{ref}
\LastPageEnding

\end{document}